\documentclass[11pt]{article}
\usepackage{latexsym, amsfonts}
\newcommand{\be}{\begin{equation}}
\newcommand{\ee}{\end{equation}}
\newcommand{\bea}{\begin{eqnarray}}
\newcommand{\eea}{\end{eqnarray}}
\newcommand{\bean}{\begin{eqnarray*}}
\newcommand{\eean}{\end{eqnarray*}}
\newcommand{\brray}{\begin{array}}
\newcommand{\erray}{\end{array}}
\newcommand{\ben}{\begin{equation}{nonumber}}
\newcommand{\een}{\end{equation}{nonumber}}

\newtheorem{dfn}{Definition}[section]
\newtheorem{thm}[dfn]{Theorem}
\newtheorem{lmma}[dfn]{Lemma}
\newtheorem{ppsn}[dfn]{Proposition}
\newtheorem{crlre}[dfn]{Corollary}
\newtheorem{xmpl}[dfn]{Example}
\newtheorem{rmrk}[dfn]{Remark}
\newcommand{\bdfn}{\begin{dfn}}
\newcommand{\bthm}{\begin{thm}}
\newcommand{\blmma}{\begin{lmma}}
\newcommand{\bppsn}{\begin{ppsn}}
\newcommand{\bcrlre}{\begin{crlre}}
\newcommand{\bxmpl}{\begin{xmpl}}
\newcommand{\brmrk}{\begin{rmrk}}
\newcommand{\edfn}{\end{dfn}}
\newcommand{\ethm}{\end{thm}}
\newcommand{\elmma}{\end{lmma}}
\newcommand{\eppsn}{\end{ppsn}}
\newcommand{\ecrlre}{\end{crlre}}
\newcommand{\exmpl}{\end{xmpl}}
\newcommand{\ermrk}{\end{rmrk}}

\newcommand{\IB}{{I\! \! B}}
\newcommand{\IC}{\mathbb{C}}

\newcommand{\IR}{\mathbb{R}}

\newcommand{\IT}{\mathbb{T}}

\newcommand{\IZ}{\mathbb{Z}}

\newcommand{\cla}{{\cal A}}
\newcommand{\clb}{{\cal B}}
\newcommand{\clc}{{\cal C}}
\newcommand{\cld}{{\cal D}}

\newcommand{\clg}{{\cal G}}
\newcommand{\clh}{{\cal H}}
\newcommand{\cli}{{\cal I}}

\newcommand{\cll}{{\cal L}}

\newcommand{\clq}{{\cal Q}}

\newcommand{\cls}{{\cal S}}

\def\a*{{\cal A}_{h,*}}
\def\B{{\cal B}(h)}
\def\B1{{\cal B}_1(h)}
\def\b{{\cal B}^{\rm s.a.}(h)}
\def\b1{{\cal B}^{\rm s.a.}_1(h)}

\newcommand{\ot}{\otimes}

\newcommand{\raro}{\rightarrow}

\newcommand{\lgl}{\langle}
\newcommand{\rgl}{\rangle}

\def \qed {$\Box$}

\begin{document}
	\[
\]
\begin{center}
{\large {\bf Quantum Isometry Groups: Examples and Computations}}\\
by\\
{\large Jyotishman Bhowmick {\footnote {The support from National Board of Higher Mathematics,  India,
 is gratefully acknowledged.}} and  Debashish Goswami{\footnote {partially supported by the project `Noncommutative Geometry and Quantum Groups' funded by the Indian National Science Academy.}}}\\
{\large Stat-Math Unit, Kolkata Centre,}\\
{\large Indian Statistical Institute}\\
{\large 203, B. T. Road, Kolkata 700 108, India}\\
{e mails: jyotish\_\ r@isical.ac.in, goswamid@isical.ac.in }\\
\end{center}
\begin{abstract}
    
    In this follow-up of \cite{goswami}, where quantum isometry group of a noncommutative manifold has been defined, we explicitly compute such quantum groups for a number of classical as well as noncommutative manifolds including the spheres and the tori. It is also proved that the    quantum isometry group of an isospectral deformation of a (classical or noncommutative) manifold is a suitable deformation of the quantum isometry group of the original (undeformed) manifold.
      
  \end{abstract} 
  \section{Introduction}
The idea of quantum isometry group of a noncommutative manifold (given by spectral triple), which has been defined by one of the authors of the present article in \cite{goswami}, is motivated by the definition and study of 
quantum permutation groups of finite sets and finite graphs by a number of mathematicians (see, e.g.
\cite{ban1}, \cite{ban2}, \cite{wang}, \cite{univ1} and references therein).
  The group of Riemannian isometries of a compact Riemannian manifold $M$ can be viewed as the universal object in the
      category of all compact metrizable groups acting on $M$, with smooth and isometric action. Therefore, to define the quantum isometry group,
       it is reasonable to  consider a
     category of compact quantum groups which act on the manifold
     (or more generally, on a noncommutative manifold given by
     spectral triple) in a `nice' way, preserving the Riemannian
     structure in some suitable sense, which is precisely formulated in \cite{goswami}, where it is also proven that a
         universal object in the category of such quantum groups does 
     exist if one makes some natural regularity assumptions on the spectral triple. Let us just sketch the definition of the quantum isometry group $\clq \equiv QISO(\cla^\infty, \clh,D)$ of a spectral triple $(\cla^\infty,\clh,D)$, without going into all the technical details, for which the reader is referred to \cite{goswami}.  The main ingredient of the definition is the Laplacian $\cll$ coming from the spectral triple (see \cite{goswami} for its construction), which coincides with  the Hodge Laplacian $-d^\ast d$ (restricted on space of smooth functions) in the classical case, where $d$ denotes the de-Rham differential. To define the Laplacian in the noncommutative case, it is assumed that
      the spectral  triple $(\cla^\infty,\clh, D)$ is of compact type 
  there is some $p>0$ such that the operator $|D|^{-p}$ (interpreted as the inverse of the restriction of $|D|^p$ on the closure of its range,
  which has a finite co-dimension since $D$ has compact resolvents) has finite nonzero Dixmier trace,
  denoted by $Tr_\omega$ (where $\omega$ is some suitable Banach limit). 
  Consider the canonical `volume form' $\tau $ coming from the Dixmier trace, i.e.
  $\tau : \clb(\clh) \raro \IC$ defined by $\tau(A):=\frac{1}{Tr_\omega(|D|^{-p})} Tr_\omega(A |D|^{-p}).$ We also assume 
  that the spectral triple
   is $QC^\infty$, i.e. ${\cla^\infty}$ and $\{ [D,a], ~a \in {\cla^\infty} \}$ are contained in the domains of all powers of the derivation $[|D|, \cdot]$. Under this assumption,
    $\tau$ is a positive faithful trace on the  $C^*$-subalgebra generated by $\cla^\infty$ and $\{ [D,a]~a \in {\cla^\infty} \}$, and using this there is a canonical construction of the Hilbert space of forms, denoted by $\clh^n_D$, $ n \geq 0$ (see \cite{fro} for details), with $\clh^0_D=L^2(\cla^\infty, \tau)$. It is assumed that the  unbounded densely defined map $d_D$ from $\clh^0_D$ to $\clh^1_D$ given by $d_D(a)=[D,a]$ for $a \in \cla^\infty$, is  closable, $\cll:=-d_D^*d_D$ has $\cla^\infty$ in its domain,  and it is left invariant by $\cll$. Moreover, we assume that $\cll$ has  compact resolvents, with its eigenvectors belonging to $\cla^\infty$, and the kernel of $\cll$ is the one-dimensional subspace spanned by the identity $1$ of $\cla^\infty$. The linear span of eigenvectors of $\cll$, which is a subspace of $\cla^\infty$, is denoted by $\cla^\infty_0$, and it is assumed that $\cla^\infty_0$ is norm-dense in the $C^*$-algebra  $\cla$ obtained by completing $\cla^\infty$.  The $\ast$-subalgebra  of $\cla^\infty$ generated by $\cla^\infty_0$ is denoted by $\cla_0$.

     It is clear that  $\cll(\cla^\infty_0) \subseteq \cla^\infty_0$,   and a compact quantum group $(\clg,\Delta)$ which has an action $\alpha$ on $\cla$ is said to act smoothly and isometrically on the noncommutative manifold $(\cla^\infty, \clh, D)$ if $({\rm id} \ot \phi) \circ \alpha(\cla^\infty_0) \subseteq \cla^\infty_0$ for all state $\phi$ on $\clg$, and also $({\rm id} \ot \phi) \circ \alpha$ commutes with $\cll$ on $\cla^\infty_{0}$. One can consider the category of all compact quantum groups acting smoothly and isometrically on $\cla$, where the morphisms are quantum group morphisms which intertwin the actions on $\cla$. It is proved in \cite{goswami}(under some regularity assumptions, which are valid for any compact connected Riemannian spin manifold with the usual Dirac operator) that there exists a universal object in this category, and this universal object is defined to be the quantum isometry group of $(\cla^\infty,\clh,D)$, denoted by $QISO(\cla^\infty, \clh, D)$, or simply as $QISO(\cla^\infty)$ or even $QISO(\cla)$ if the spectral triple is understood. In fact, we have considered a bigger category, namely the category of `quantum families of smooth isometries' (see \cite{goswami} for details), which is motivated by the ideas of Woronowicz and Soltan (\cite{woro_pseudo}, \cite{soltan}), and identified the underlying $C^*$-algebra of the quantum isometry group as a universal object in this bigger category.
     
        We believe that a detailed study of quantum
      isometry groups will not only give many new and interesting
      examples of compact quantum groups, it will also contribute
      to the understanding of quantum group covariant spectral
      triples. For this, it is important to explicitly describe quantum isometry groups of sufficiently many classical and noncommutative manifolds. This is our aim in this paper. We have computed quantum isometry groups of classical and noncommutative spheres and tori, and also obtained a gereral principle for computing such quantum groups, by proving that the quantum isometry group of an isospectral deformation of a (classical or noncommutative) manifold is a deformation of the quantum isometry group of the original (undeformed) manifold.
      
 Throughout the paper, we have denoted by $\cla_1 \ot \cla_2$ the minimal (injective) $C^*$-tensor product between two $C^*$-algebras $\cla_1$ and $\cla_2$. The symbol $\ot_{\rm alg}$ has been used to denote the algebraic tensor product between vector spaces or algebras.  
 
 For a compact quantum group $\clg$, the dense unital $\ast$-subalgebra generated by the matrix coefficients of irreducible  unitary representations has been denoted by $\clg_0$. The coproduct of $\clg$, say $\Delta$, maps $\clg_0$ into the algebraic tensor product $\clg_0 \ot_{\rm alg} \clg_0$, and there exist canonical antipode and counit defined on $\clg_0$ which make it into a Hopf $\ast$-algebra ( see \cite{woro} for the details ).      
  
  \section{Computation of the quantum isometry groups of the sphere and tori}
  \subsection{Computation for the commutative spheres}
  Let $\clq$ be the quantum isometry group of $S^2$ and let $\alpha$ be the action of $\clq$ on $C(S^2)$.
      Let $\cll$ be the Laplacian on $S^2$ given by $$\cll=\frac{\partial^2}{\partial \theta^2}+
{\rm cot}(\theta) \frac{\partial}{\partial \theta}+\frac{1}{{\rm sin}^2(\theta)}\frac{\partial^2}{\partial \psi^2},$$ and   the cartesian coordinates $x_1$, $x_2$, $x_3$ for $S^2$ are given by $x_1=r\cos{\psi}  \sin{\theta}$, $x_2=r\sin{\psi} \sin{\theta}$, $x_3=r\cos{\theta}$. In the cartesian coordinates, $\cll=\sum_{i=1}^3 \frac{\partial^2}{\partial x_i^2}.$
      
      The eigenspaces of $\cll$ on $S^2$ are of the form $$E_k={\rm Sp}\{(c_1X_1+c_2X_2+c_3X_3)^k~:~c_i\in \IC,i=1,2,3,~ \sum c_i^2=0\},$$ where $k \geq 1$. $E_k$ consists of harmonic homogeneous polynomials of degree $k$ on $R^3$ restricted to $S^2$.( See \cite{Helgason}, page 29-30 ).
      
      We begin with the following lemma, which says that any smooth isometric action by a quantum group must be `linear'. 
      
      \blmma
      The action  $\alpha$ satifies $\alpha(x_i)=\sum_{j=1}^{3} x_j\otimes Q_{ij}$ where $Q_{ij} \in \clq, i=1,2,3$.
      \elmma
    {\it Proof :}\\  
      $\alpha$ is a smooth isometric action of $\clq$ on $C(S^2)$, so $\alpha$ has to preserve the eigenspaces of the laplacian $\cll$. In particular, it has to preserve $E_1={\rm Sp}\{ c_1x_1+c_2x_2+c_3x_3~:~c_i \in \IC,i=1,2,3,      \sum_{i=1}^{3}c^2_i=0\}.$
   
   Now note that $x_1+ix_2,  x_1-ix_2 \in E_1$, hence $x_1,x_2\in E_1.$
          Similarly $x_3 \in E_1$ too. 
     Therefore $E_1={\rm Sp}\{ x_1,x_2,x_3 \}$, which completes the proof of the lemma.
     \qed
     
     Now, we state and prove the main result of this section, which identifies $\clq$ with the commutative $C^*$ algebra of continuous functions on the isometry group of $S^2$, i.e. $O(3)$.
     \bthm
     The quantum isometry group $\clq$ is commutative as a $C^*$ algebra.
     \ethm
     {\it Proof :}\\          
     We begin with the expression $$ \alpha(x_i)=\sum_{j=1}^3 x_j \otimes Q_{ij},~i=1,2,3,$$ and also note that $x_1,x_2,x_3$ form a basis of $E_1$ and $\{ x_1^2,x_2^2,x_3^2,x_1x_2,x_1x_3, x_2x_3 \}$ is a basis of $E_2$. 
          Since $x_i^*=x_i$ for each $i$ and $\alpha$ is a $\ast$-homomorphism, we must have $Q_{ij}^*=Q_{ij} ~ \forall i,j=1,2,3$. 
           Moreover, the condition $x^2_1+x^2_2+x^2_3=1$ and the fact that $\alpha$ is a homomorphism gives:
           
                     $$ Q^2_{1j}+Q^2_{2j}+Q^2_{3j}=1,~\forall j=1,2,3.$$

          Again,the condition that $x_i$,$x_j$ commutes $\forall i,j$ gives \be \label{2c} Q_{ij}Q_{kj}=Q_{kj}Q_{ij} \forall i,j,k, \ee
          \be \label{3c} Q_{ik}Q_{jl}+Q_{il}Q_{jk}=Q_{jk}Q_{il}+Q_{jl}Q_{ik}.\ee

               Now, it follows from the  Lemma 2.12 in \cite{goswami} 
                that $\tilde{\alpha}: C(S^2) \otimes \clq \raro C(S^2) \otimes \clq$ defined by $\tilde{\alpha}(X \otimes Y)=\alpha(X)(1 \otimes Y)$ extends to  a unitary of the Hilbert $\clq$-module $L^2 ( S^2 ) \otimes \clq$ (or in other words, $\alpha$ extends to a unitary representation of $\clq$ on $L^2(S^2)$). 
               But  $\alpha$ keeps $V={\rm Sp}\{x_1,x_2,x_3\}$ invariant.
               So $\alpha$ is a unitary representation  of $\clq$ on $V$, i.e. $Q = (( Q_{ij} ))  \in M_3 ( \clq )$ is a unitary,      hence $Q^{-1}=Q^*=Q^T$, since in this case entries of $Q$ are self-adjoint elements.

               Clearly, the matrix $Q$ is a $3$-dimensional unitary representation of $\clq$. Recall that ( cf \cite{VanDaele} ) 
            the antipode $\kappa$ on the matrix elements of a finite-dimensional   unitary representation $U^\alpha \equiv ( u_{pq}^\alpha)$ is given by $\kappa (u_{pq}^\alpha ) =( u_{qp}^\alpha )^* .$
               
               So we obtain \be \label{5c} \kappa( Q_{ij} )= Q^{-1}_{ij}=Q^T_{ij}=Q_{ji}.\ee
               
               Now from ( \ref{2c} )  , we have $Q_{ij}Q_{kj} = Q_{kj}Q_{ij}.$
               Applying $ \kappa$ on this equation and using the fact that $\kappa$ is an antihomomorphism along with ( \ref{5c} ) , 
            we have $Q_{jk}Q_{ji} = Q_{ji}Q_{jk}$
               Similarly , applying $\kappa$ on ( \ref{3c} ), we get 
               
               $$ Q_{lj}Q_{ki} + Q_{kj}Q_{li} = Q_{li}Q_{kj} +Q_{ki}Q_{lj}~ \forall i,j,k,l.$$
         Interchanging between $k$ and $i$ and also between $l,j$ gives 
               \be \label{6c} Q_{jl}Q_{ik} +Q_{il}Q_{jk} =Q_{jk}Q_{il} +Q_{ik}Q_{jl}~ \forall i,j,k,l.\ee
               
               Now, by (\ref{3c} )-( \ref{6c} ) , we have
               
               $$ [ Q_{ik},Q_{jl} ] =[ Q_{jl},Q_{ik} ],$$ 
               
               hence $$ [ Q_{ik},Q_{jl} ] = 0.$$
               
                            Therefore the entries of  the matrix $Q$ commute among themselves. 
            However, by faithfulness of the action of $\clq$, it is clear that the $C^*$-subalgebra generated by entries of $Q$ (which forms a quantum subgroup of $\clq$ acting on $C(S^2)$ isometrically) must be the same as  $\clq$,  so $\clq$ is commutative.
            \qed
               
               So $\clq=C( G )$ for some compact group $G$ acting by isometry on $C (S^2 )$ and is also universal in this                   category, i.e. $\clq=C( O( 3 )).$

     \brmrk
      
        Similarly, it can be shown that $ QISO ( S^n ) $ is commutative for all $ n \geq 2. $
        
        \ermrk

      \subsection {The commutative one-torus}
       
      Let $\clc=C(S^1)$ be the $C^*$-algebra of continuous functions on the one-torus $S^1$. Let us denote by $z$ and $\overline{z}$ the identity function (which is the generator of $C(S^1)$) and its conjugate respectively. The Laplacian coming from the standard Riemannian metric is given by $\cll(z^n)=-n^2 z^n$, for $n \in \IZ$, hence the eigenspace corresponding to the eigenvalue $-1$ is spanned by $z$ and $\overline{z}$ . Thus, the action of a compact quantum group acting smoothly and isometrically (and faithfully) on $C(S^1)$ must be {\it linear} in the sense that its action must map $z$ into an element of the form $z \ot A +\overline{z} \ot B$. However, we show below that this forces the quantum group to be commutative as a $C^*$ algebra, i.e. it must be the function algebra of some compact group .
       
       \bthm
       Let $\alpha $ be a faithful, smooth and  linear action of a compact quantum group $(\clq,\Delta)$ on $C ( S^1 )$ defined by $ \alpha ( z ) = z \otimes A + \overline{z} \otimes B$. 
      Then $\clq$ is a commutative $C^*$ algebra.
      \ethm 
      {\it Proof :}\\
    By the assumption of faithfulness, it is clear that $\clq$ is generated (as a unital $C^*$ algebra) by $A$ and $B$. Moreover, recall that smoothness in particular means that $A$ and $B$ must belong to the algebra $\clq_0$ spanned by matrix elements of irreducible representations of $\clq$ . Since $z \overline{z} =\overline{z}z =1$ and $\alpha$ is a $\ast$-homomorphism, we have 
        $ \alpha ( z ) \alpha (  \overline{z} ) = \alpha ( \overline{z} ) \alpha ( z ) = 1 \otimes 1 $.

    Comparing coefficients of $z^2,{ \overline {z}}^2$ and $1$ in both hand sides of the relation  $ \alpha ( z ) \alpha (  \overline{z} ) =1 \otimes 1$, we get
    
    \be \label{1d} 
    AB^* = BA^* = 0,~~~
         AA^* + BB^* =1.\ee
    
    Similarly, $ \alpha (\overline{z} ) \alpha ( z ) = 1 \otimes 1$ gives 
    
    \be \label{2d} B^*A = A^*B =0,~~~
         A^*A +B^*B =1. \ee
    
    Let $U =A+B$ ,
           $P=A^*A$ ,         
        $Q=AA^*$.      
        Then it follows from  (\ref{1d}) and (\ref{2d}) that $U$ is a unitary  and $P$ is a projection since  $P$ is self adjoint and  \bean \lefteqn{ P^2} ~~
        & =& A^*AA^*A ~~       
           = A^*A( 1-B^*B )~~
           = A^*A - A^*AB^*B ~~
            = A^*A  ~~                           
           =P.\eean 
                                        
          Moreover ,  \bean \lefteqn{ UP}\\
          & =& ( A + B ) A^*A~~
                             =AA^*A + BA^*A ~~                   
                   =AA^*A \\&& (~{\rm since}~ BA^* =0~{\rm from}~(\ref{1d}) )\\ 
                   &=& A ( 1-B^*B )
                    ~~=A-AB^*B ~~=A .\eean             
       Thus, $A=UP$ , $B=U-UP=U(1-P)\equiv UP^\perp$ , so $\clq=C^*(A,B)=C^*(U,P)$.

  We can rewrite the action $\alpha$ as follows:\\
    
    $$ \alpha ( z ) =z \otimes UP + \overline {z} \otimes UP^\bot.$$
    
   The coproduct $\Delta$ can easily be calculated from the requirement 
     $ (id \otimes \Delta )\alpha =( \alpha \otimes id ) \alpha$ , and it is given by :
    
    
        
    
   \be \label{6d} \Delta ( UP ) = UP \otimes UP +P^ \bot U^{-1 } \otimes UP^ \bot  \ee   
    \be \label{7d} \Delta ( UP^ \bot ) = UP^\bot \otimes UP + PU^{-1} \otimes UP^ \bot. \ee   
  From this, we get
   \be \label{8d}  \Delta ( U ) =U \otimes UP +U^{-1} \otimes UP^ \bot, \ee   
 \be \label{9d} \Delta ( P ) = \Delta (U^{ -1 }) \Delta ( UP )
                          =P \otimes P +UP^\bot U^{-1} \otimes P^\bot.\ee 
                       
   It can be checked that $ \Delta $ given by  the above expression is coassociative.
   
   Let $h$ denote the right-invariant Haar state on $\clq$. By the general theory of compact quantum groups, $h$ must be faithful on $\clq_0$.  
      We have (by right-invariance of $h$): $$ ({\rm id} \otimes h) ( P \otimes P + UP^\bot U^{-1} \otimes P^\bot ) =h( P )1.$$
   
  That is, we have \be \label{10d} h( P^\bot )UP^ \bot U^{-1} = h ( P )P^ \bot.\ee
   
   Since $P$ is a positive element in $\clq_0$ and $h$ is faithful on $\clq_0$, $h(P)=0$ if and only if $P=0$. Similarly , $h(P^\bot)=0$, i.e. $h(P)=1$, if and only if $P=1$. However, if $P$ is either $0$ or $1$, clearly  $\clq =C^*( U,P )=C^*( U )$, which is commutative. On the other hand, if we assume that $P$ is not a trivial projection, then $h(P)$ is strictly between $0$ and $1$, and we have 
 from  ( \ref{10d} ) $$  UP^\bot U^{-1} = \frac{h( P )}{ 1-h( P )} P^\bot .$$
   
   Since both $UP^ \bot U^{-1} $ and $P^ \bot $ are nontrivial projections, they can be scalar multiples of each other if and only if they are equal, so we conclude that 
   $ UP^\bot U^{-1}=P^\bot$, i.e. $U$ commutes with $P^\bot$, hence with $P$, and $\clq$ is commutative.
  \qed

   \subsection{Commutative and noncommutative two-tori} 
     Fix a real number $\theta$, and let $ \cla_ \theta $ be the universal $ C^{*} $ algebra generated by two unitaries $ U $ and $ V $ such that $ U V = \lambda V U $, where $\lambda:=e^{2 \pi i \theta}$. It is well-known (see \cite{con}) that       the set $ \{ U^{m}V^{n} : m,n \in \IZ \}  $ is  an orthonormal basis for  $ L^{2} ( \cla_{ \theta } , \tau ),     $ where $ \tau$ denotes the unique faithful normalized trace on $\cla_{\theta}$ given by, $\tau ( \sum a_{m n} U^{m} V^{n} ) = a_{0 0} $.     
      We shall denote by  $ \left\langle A , B \right\rangle = \tau ( A^{*} B )$    the inner product on $\clh_0:=L^2(\cla_\theta,\tau)$. Let $\cla_\theta^{\rm fin}$ be the unital  $\ast$-subalgebra generated by finite complex linear combinations of $U^mV^n$, $m,n \in \IZ$, and $d_1,d_2$ be the maps on $\cla^{\rm fin}_\theta$ defined by $d_1(U^mV^n)=mU^mV^n$, $d_2(U^mV^n)=nU^mV^n)$. We consider the canonical spectral triple (see \cite{con} for details) $(\cla^{\rm fin}_\theta, \clh, D)$, where $\clh=\clh_0 \oplus \clh_0$, $D=\left( \begin{array}{cc} 0 & d_1+id_2 \\ d_1-id_2 & 0 \end{array} \right),$ and the representation of $\cla_\theta$ on $\clh$ is the diagonal one,   i.e. $a \mapsto \left( \begin{array}{cc} a & 0 \\ 0 & a \end{array} \right).$ Clearly, the corresponding Laplacian $\cll$ is given by $\cll(U^mV^n)=-(m^2+n^2) U^mV^n,$ and it is also easy to see that the algebraic span of eigenvectors of $\cll$ is nothing but the space $\cla^{\rm fin}_\theta$, and moreover, all the assumptions in \cite{goswami} required for defining the quantum isometry group are satisfied. 
         
         Let $\clq$ be the quantum isometry group of the above spectral triple, with the 
  smooth isometric action of  on $ \cla_ \theta $ given by $\alpha : \cla_\theta \raro \cla_\theta \ot \clq$. 
     By definition,  $ \alpha $  must keep invariant the eigenspace of $ \cll $ corresponding to the eigen value $- 1 $ , spanned by $ U,V,U^{-1},V^{-1} $.Thus, the action $ \alpha $ is given by:   
    $$  \alpha ( U ) = U \otimes A_{1} + V \otimes B_{1} + U^{-1} \otimes C_{1} + V^{-1} \otimes D_{1},$$
  $$ 
   \alpha ( V ) = U \otimes A_{2} + V \otimes B_{2} + U^{-1} \otimes C_{2} + V^{-1} \otimes D_{2},$$ 
   for some $ A_{i},B_{i},C_{i},D_{i} \in \clq ,i=1,2 $, and by faithfulness of the action of quantum isometry group (see \cite{goswami}), the norm-closure of the unital $\ast$-algebra generated by $A_i,B_i,C_i,D_i;i=1,2$ must be the whole of $\clq$.
   
   Next we derive a number of conditions on $ A_{i},B_{i},C_{i},D_{i} , i=1,2 $ using the fact that $ \alpha $ is a $ \ast $   homomorphism.
  \blmma    
    \label{Lemma 1}

   The condition $ U^{*} U = 1 = U U^{*} $ gives:
   
  \be \label{lem1.1a}  A^{*}_{1} A_{1} +  B^{*}_{1}B _{1} +  C^{*}_{1}C_{1} + D^{*}_{1}D_{1} = 1    \ee

  \be \label{lem1.2} A^{*}_{1} B_{1} + \lambda D^{*}_{1} C_{1} = A^{*}_{1} D_{1} + \overline { \lambda } B^{*}_{1} C_{1} = 0 \ee

  \be \label{lem1.3} C^{*}_{1} D_{1} + \lambda B^{*}_{1} A_{1} = C^{*}_{1} B_{1} + \overline { \lambda } D^{*}_{1} A_{1} = 0 \ee
   
  \be \label{lem1.4} A^{*}_{1} C_{1} = B^{*}_{1} D_{1} = C^{*}_{1} A_{1} = D^{*}_{1} B_{1} = 0 \ee

  \be \label{lem1.5} A_{1} A^{*}_{1} + B_{1} B^{*}_{1} + C_{1}C^{*}_{1} + D_{1}D^{*}_{1} = 1 \ee
   
   \be \label{lem1.6}A_{1}B^{*}_{1} + \lambda D_{1}C^{*}_{1} = A_{1} D^{*}_{1} + \overline{ \lambda } B_{1}C^{*}_{1} = 0  \ee

  \be \label{lem1.7}C_{1} D^{*}_{1} + \lambda B_{1}A^{*}_{1} = C_{1}B^{*}_{1} + \overline { \lambda } D_{1}A^{*}_{1} = 0 \ee
   
  \be \label{lem1.8} A_{1}C^{*}_{1} = B_{1}D^{*}_{1} = C_{1}A^{*}_{1} = D_{1}B^{*}_{1} = 0 \ee

   \elmma
 {\it Proof :}\\
  We get ( \ref{lem1.1a} ) - ( \ref{lem1.4} ) by using the condition $ U^{*} U = 1 $ along with the fact that $ \alpha $ is a homomorphism and then comparing the coefficients of $ 1, U^{*}V ,{ U^{*} }^{2}, U^{*}V^{*} , U V^{*} ,{ V^{*} }^{2} , U^{2} , U V , V^{2} . $
   
   Similarly the condition $ U U^{*} = 1 $ gives ( \ref{lem1.5} )  -( \ref{lem1.8} ).\qed

    \blmma    
    \label{Lemma 2}

    We have analogues of ( \ref{lem1.1a} ) - ( \ref{lem1.8} ) with $ A_{1},B_{1},C_{1},D_{1} $  replaced by $ A_{2},B_{2},C_{2},D_{2} $ respectively.
  
   \elmma
  {\it Proof :}\\ We use the condition $ V^{*}V = V V^{*} = 1 $ \qed \\
  
   Now we note that if $ \alpha ( U^{m} V^{n} ) = \sum c_{kl}  U^{k} V^{l} \otimes Q_{kl} $ for some $ Q_{kl} \in \clq $, then the condition that $ \alpha $ commutes with the laplacian implies    
   $ c_{kl} = 0 $ unless $ k^{2} + l^{2} = m^{2} + n^{2} $.
   
   We use this observation in the next lemma.
   
   \blmma
   \label{Lemma 3}
   
   Inspecting the terms with zero coefficient in $ \alpha ( U^{*} V ) , \alpha ( V U^{*} ) , \alpha ( U V ) , \alpha ( V U ) $, we get 
   
   \be \label{lem3.1}  C^{*}_{1} A_{2} = 0 , D^{*}_{1} B_{2} = 0 , A^{*}_{1} C_{2} = 0 , B^{*}_{1}D_{2} = 0 \ee

   \be \label{lem3.2}  A_{2}C^{*}_{1} = 0 , B_{2}D^{*}_{1} = 0 , C_{2} A^{*}_{1} = 0 , D_{2} B^{*}_{1} = 0 \ee

   \be \label{lem3.3}  A_{1} A_{2} = 0 , B_{1} B_{2} = 0 , C_{1}C_{2} = 0 , D_{1} D_{2} = 0 \ee

    \be \label{lem3.4}  A_{2} A_{1} = 0 ,B_{2} B_{1} = 0, C_{2} C_{1} = 0, D_{2} D_{1} = 0.  \ee 
   \elmma
   
    {\it Proof :}\\   
   The equation ( \ref{lem3.1} ) is  obtained from the coefficients of $  U^{2}, V^{2} ,{ U^{*} }^{2} , { V^{*} }^{2} $  in $ \alpha ( U^{*} V ) $ while ( \ref{lem3.2} ) , ( \ref{lem3.3} ) , ( \ref{lem3.4} ) are obtained from the same coefficients in $ \alpha ( V U^{*} ) , \alpha ( U V ) , \alpha ( V U ) $ respectively. \qed
   
   \blmma
   \label{Lemma 4}:

  $  A_{1} B_{2} + \overline{ \lambda } B_{1} A_{2}  = \lambda A_{2} B_{1} + B_{2} A_{1}  $
   
 $ A_{1} D_{2} +  \lambda D_{1} A_{2}  = \lambda A_{2} D_{1} + { \lambda }^2 D_{2} A_{1} $
   
 $ C_{1} B_{2} +  \lambda B_{1} C_{2} = \lambda C_{2} B_{1} + { \lambda }^{2} B_{2} C_{1} $
   
 $ C_{1} D_{2} + \overline{ \lambda } D_{1} C_{2}  = \lambda C_{2} D_{1} +  D_{2} C_{1} $
   
   \elmma
   
   {\it Proof :}\\ Follows from the relation $ \alpha ( U V ) = \lambda \alpha ( V U ) $ and equating non zero coefficients of $ U V , U V^{-1} , U^{-1} V $ and $ U^{-1} V^{-1} $. \qed \vspace{4mm}

      Now, by Lemma 2.12 in \cite{goswami} it follows that
            $ \tilde{\alpha}: \cla_{\theta} \otimes \clq \raro \cla_{ \theta } \otimes \clq $ defined by                $\tilde{\alpha}(X \otimes Y)=\alpha(X)(1 \otimes Y)$ extends to  a unitary of the Hilbert $\clq$-module $L^2 (             \cla_{\theta},\tau ) \otimes \clq$ (or in other words, $\alpha$ extends to a unitary representation of $\clq$ on                 $L^2(\cla_{\theta},\tau)$). 
               But $\alpha$ keeps $W = {\rm Sp}\{U,V,U^{*},V^{*}\}$ invariant( as observed in the beginning of this section ).
               So $\alpha$ is a unitary representation of $\clq$ on $W$.Hence, the matrix ( say M ) corresponding to the 4 dimensional representation of $ \clq $ on $W$ is a unitary in $ M_4 ( \clq )$.
               
      From the definition of the action it follows that $ M =\left(  \begin {array} {cccc}
       A_{1} &  A_{2} & C^{*}_{1} &  C^{*}_{2} \\ B_{1} & B_{2} & D^{*}_{1} & D^{*}_{2} \\ C_{1} & C_{2} & A^{*}_{1} & A^{*}_{2} \\ D_{1} & D_{2} & B^{*}_{1} & B^{*}_{2} \end {array} \right )   $
      
      Since $ M $ is the matrix corresponding to a finite dimensional unitary representation, $ \kappa ( M_{k l } )= M^{-1}_{ k l } $ where $ \kappa $ denotes the antipode of $ \clq $ (See \cite{VanDaele})
      
    But $ M $ is a unitary, $ M^{-1} = M^{*} $
    
  So,$  ( k ( M_{k l} ) ) = \left ( \begin {array} {cccc}
   A^{*}_{1} & B^{*}_{1} & C^{*}_{1} & D^{*}_{1} \\ A^{*}_{2} & B^{*}_{2} & C^{*}_{2} & D^{*}_{2} \\ C_{1} & D_{1} & A_{1} & B_{1} \\ C_{2} & D_{2} & A_{2} & B_{2} \end {array} \right ) $

       \blmma
   \label{Lemma 5a}:

      $ A_{1} $ is a normal partial isometry and hence has same domain and range.
      \elmma
      
   {\it Proof :}\\ From the relation $ A^{*}_{1} A_{1} +  B^{*}_{1}B _{1} + C^{*}_{1}C_{1} + D^{*}_{1}D_{1} = 1 $ in Lemma \ref{Lemma 1}, we have by applying $ \kappa $, $ A^{*}_{1} A_{1} +  A^{*}_{2}A _{2} + C_{1}C^{*}_{1} + C_{2}C^{*}_{2} = 1 $ .       
      Applying $ A_{1} $ on the right of this equation and using $ C^{*}_{1} A_{1} = 0 $ from Lemma  \ref{Lemma 1}, and $ A_{2}A_{1} = A^{*}_{1} C_{2} = 0 $ from Lemma \ref{Lemma 3}, we have 
      
  \be \label {i}    A^{*}_{1} A_{1} A_{1} = A_{1}  \ee 
     
     Again, from the relation $ A_{1} A^{*}_{1} +  B_{1}B^{*}_{1} + C_{1}C^{*}_{1} + D_{1}D^{*}_{1} = 1 $ in Lemma  \ref{Lemma 1}, applying $ \kappa $ and multiplying by $ A^{*}_{1} $ on the right, and then using $ C_{1} A^{*}_{1} = 0 $ from Lemma  \ref{Lemma 1},     
    $  A_{1} A_{2} = C_{2} A^{*}_{1} = 0 $ from Lemma \ref{Lemma 3} , we have \be \label {ii}  A_{1} A^{*}_{1}A^{*}_{1} = A^{\ast}_{1}  \ee
     
     From (\ref {i}), we have \be \label {iii}  ( A^{*}_{1} A_{1} ) ( A_{1} A^{*}_{1} ) = A_{1} A^{*}_{1}  \ee 
     
     By taking $*$ on (\ref {ii}), we have \be \label {iv}  A_{1} A_{1} A^{*}_{1} = A_{1}  \ee
     
     So, by multiplying by $ A^{*}_{1} $ on the left,  we have
     
    \be \label {v} ( A^{*}_{1} A_{1} )( A_{1} A^{*}_{1} ) = A^{*}_{1} A_{1}   \ee
      
      From (\ref {iii})  and (\ref {v}), we have      
     $ A_{1} A^{*}_{1} = A^{*}_{1} A_{1} $, i.e $ A_{1} $ is normal.     
     So, $ A_{1} = A^{*}_{1} A_{1} A_{1} $ ( from ( \ref {i} ) )     
                $ = A_{1} A^{*}_{1} A_{1} $
                
      Therefore, $ A_{1} $ is a partial isometry which is normal and hence has same domain and range.\qed
      
      \brmrk

      In an exactly similar way, it can be proved that $ D_{1} $ is a normal partial isometry and hence has same domain and range. 
      \ermrk
      
    \blmma
   \label{Lemma 6}:

   We have $ C^{*}_{1} B^{*}_{1} = C^{*}_{2} B^{*}_{2} = A_{1} D_{1} = A_{2} D_{2} = B_{1} C^{*}_{1} =   B^{*}_{1} C^{*}_{1} = B_{1} A_{1} = A_{1} B^{*}_{1} = D_{1} A_{1} = A^{*}_{1} D_{1} = C^{*}_{1} B_{1} = D_{1} C^{*}_{1} = 0 $
  \elmma
  
 {\it Proof :}\\ Using $ A_{2} C^{*}_{1} = B_{2} D^{*}_{1} = C_{2} A^{*}_{1} = D_{2} B^{*}_{1} = 0 $  from Lemma \ref{Lemma 3} and applying $ \kappa $, we have the first four equalities. 
 But $ A_{1} D_{1} = 0 $ .Hence \be \label{lem6.1} Ran ( D_{1} ) \subseteq Ker A_{1}  \ee  
 
 By the above made remark, \be \label{lem6.2}   Ran ( D_{1} ) = Ran ( D^{*}_{1} )   \ee
 
 ( \ref{lem6.1} )  and ( \ref{lem6.2} ) imply  $ Ran ( D^{*}_{1} ) \subseteq Ker ( A_{1} ) $, so 
                           $ A_{1} D^{*}_{1} = 0 $
                          
   But from Lemma \ref{Lemma 1}, we have $ A_{1} D^{*}_{1} + \overline{ \lambda } B_{1} C^{*}_{1} = 0 $, which gives $ B_{1} C^{*}_{1} = 0 $.
   
   From Lemma \ref{Lemma 3}, we have $ C^{*}_{1} A_{2} =  A_{2} A_{1} = 0 $, from which it follows by applying $ \kappa $, that $ B^{*}_{1} C^{*}_{1} = A^{*}_{1} B^{*}_{1} = 0 $     
So, $ B_{1} A_{1} = 0 $.
                         
    So, $ Ran ( A_{1} ) \subseteq Ker B_{1} $
    
    But by Lemma \ref{Lemma 5a}, $ A_{1} $ is a normal partial isometry and so has same range and domain.
    
    Thus, $ Ran ( A^{*}_{1} ) \subseteq Ker ( B_{1} ) $ which implies $ B_{1} A^{*}_{1} = 0 $ , i.e ,   
    \be \label{lem6.3}A_{1} B^{*}_{1} = 0  \ee    
      
    Again, from Lemma \ref{Lemma 3}, $ A^{*}_{1} C_{2} = 0 $. Hence, by applying $ \kappa $, $ D_{1} A_{1} = 0 $ , i.e, $ A^{*}_{1} D^{*}_{1} = 0 $.    
    But $ D_{1} $ is a partial isometry ( from the remark following Lemma \ref{Lemma 5a} ), we conclude  $ A^{*}_{1} D_{1} = 0 $.    
    But by Lemma \ref{Lemma 1}, we have      
    $ A^{*}_{1} D_{1} + \overline { \lambda } B^{*}_{1} C_{1} = 0 $    
    But, $A^{*}_{1} D_{1} = 0 $ implies $ B^{*}_{1} C_{1} = 0 $ i.e, $ C^{*}_{1} B_{1} = 0 $    
    Also, $ A_{1} B^{*}_{1} = 0 $ ( from (\ref{lem6.3} ) ) and $ A_{1} B^{*}_{1} + \lambda D_{1} C^{*}_{1} = 0 $( by Lemma \ref{Lemma 1} ), so 
   $  D_{1} C^{*}_{1} = 0 $
    \qed
   
   \blmma
   \label{Lemma 7}:

   $ C_{1} $ is a normal partial isometry and hence has same domain and range.
   \elmma
   
   {\it Proof :}\\From the relation   
         $ A^{*}_{1} A_{1} +  B^{*}_{1}B _{1} + C^{*}_{1}C_{1} + D^{*}_{1}D_{1} = 1 $ in Lemma \ref{Lemma 1}, multiplying by $ C^{*}_{1} $ on the right and using $ A_{1} C^{*}_{1} = 0 $ from Lemma \ref{Lemma 1}, and $ B_{1}C^{*}_{1} = D_{1} C^{*}_{1} = 0 $ from Lemma \ref{Lemma 6}, we have \be \label{lem7.1}  C^{*}_{1} C_{1} C^{*}_{1} = C^{*}_{1} \ee
         
 Therefore, $ C^{*}_{1} $ and hence $ C_{1} $ is a  partial isometry.

 Also,from Lemma \ref{Lemma 1}, $  A^{*}_{1} A_{1} +  B^{*}_{1}B _{1} + C^{*}_{1}C_{1} + D^{*}_{1}D_{1} = 1 =  A_{1} A^{*}_{1} + B_{1} B^{*}_{1} + C_{1}C^{*}_{1} + D_{1}D^{*}_{1} $ .       
       Using the normality of $ A_{1} $ and $ D_{1} $ ( obtained from Lemma \ref{Lemma 5a} and the remark following it ) to this equation, we have \be \label{lem7.2} B^{*}_{1}B _{1} + C^{*}_{1}C_{1} =  B_{1} B^{*}_{1} + C_{1}C^{*}_{1} \ee   
       
       Multiplying by $ C^{*}_{1} $ to the left of (\ref{lem7.2}), and using $ C^{*}_{1} B^{*}_{1} = C^{*}_{1} B_{1} = 0 $ from Lemma \ref{Lemma 6}, we have :      
                   $ C^{*}_{1} C^{*}_{1} C_{1} = C^{*}_{1} C_{1} C^{*}_{1}. $                   
                   But $ C^{*}_{1} C_{1} C^{*}_{1} = C^{*}_{1} $ ( from (\ref{lem7.1}) )                   
                   , hence $ C^{*}_{1} C^{*}_{1} C_{1} = C^{*}_{1} $                   
   Applying $ C_{1} $ on the left, we have \be \label{lem7.3} ( C_{1} C^{*}_{1} ) ( C^{*}_{1} C_{1} ) = C_{1} C^{*}_{1}  \ee 
   
   Now multiplying by $ C^{*}_{1} $ on the right of (\ref{lem7.2}) and using $ B_{1} C^{*}_{1} = B^{*}_{1} C^{*}_{1} = 0 $ from Lemma \ref{Lemma 6}, we have $ C^{*}_{1} C_{1} C^{*}_{1} = C_{1} C^{*}_{1} C^{*}_{1} $   
        and using (\ref{lem7.1}), we have $ C_{1} C^{*}_{1} C^{*}_{1} = C^{*}_{1} . $        
        Thus,  $ C^{*}_{1} C_{1} =( C_{1} C^{*}_{1}) ( C^{*}_{1} C_{1} )         
                              = C_{1} C^{*}_{1} $ ( by (\ref{lem7.3}) ),                               
        hence $ C_{1} $ is a normal partial isometry and so has the same domain and range. \qed

      \brmrk :
      
     1. In the same way, it can be proved that $ B_{1} $ is a normal partial isometry and hence has same domain and range.
     
     2. In an exactly similar way, it can be proved that $ A_{2} , B_{2} , C_{2} , D_{2} $ are  normal partial isometries and hence has same domain and range.
     \ermrk
     
     \blmma
   \label{Lemma 8}:   
     
     We have   $ A_{1} C_{2}= B_{1} D_{2} = C_{1} A_{2} = D_{1} B_{2} . $
   \elmma

    {\it Proof :}\\           
  By Lemma \ref{Lemma 6}, we have $ A_{1} D_{1} = A_{2} D_{2} = C^{*}_{2} B^{*}_{2} =  0 $.  
  Now, using the fact that $ D_{1}, D_{2} $ and $ B_{2} $ are normal partial isometries, we have $ A_{1} D^{*}_{1} = A_{2} D^{*}_{2} = C^{*}_{2} B_{2} = 0 $  
  Taking adjoint and applying $ \kappa $, we have the first, second and the fourth equalities.  
  To prove the third one, we take adjoint of the relation $ C^{*}_{1} B_{1} = 0 $  obtained from Lemma \ref{Lemma 6} and then apply $ \kappa $ . \qed \vspace{4mm}

    Now we define for $ i = 1,2 $
    
    $ A^{*}_{i} A_{i} = P_{i}, B^{*}_{i} B_{i} = Q_{i}, C^{*}_{i} C_{i} = R_{i}, S_{i} = 1 - P_{i} - Q_{i} - R_{i} $
     
    $  A_{i} A^{*}_{i} = P^{{\prime}}_{i}, B_{i} B^{*}_{i} = Q^{{\prime}}_{i}, C_{i} C^{*}_{i} = R^{{\prime}}_{i}, S^{{\prime}}_{i} = 1 - P^{{\prime}}_{i} - Q^{{\prime}}_{i} - R^{{\prime}}_{i} $
           
     By Lemma \ref{Lemma 1}, and the remark following it, we have $ D^{*}_{i} D_{i} = 1  - ( P_{i} + Q_{i} + R_{i} ) $   
                                                     and $ D_{i} D^{*}_{i} = 1 - ( P^{{\prime}}_{i} + Q^{{\prime}}_{i} + R^{{\prime}}_{i} ) $
      
      Also we note that, since $ A_{i}, B_{i}, C_{i}, D_{i} $ are normal, it follows that $ P_{i}= P^{{\prime}}_{i}, Q_{i} = Q^{{\prime}}_{i}, R_{i} = R^{{\prime}}_{i}, S_{i} =  S^{{\prime}}_{i} $.
                                                    
      \blmma :
      \label{Lemma 9} 
     $ P_{1} + R_{1} = 1 - ( P^{{\prime}}_{2} + R^{{\prime}}_{2} ) $
      \elmma 
      
     {\it Proof :}\\ 
                 From Lemma \ref{Lemma 3}, $ A_{1} A_{2} = B_{1} B_{2} = C_{1} C_{2} = D_{1} D_{2} = 0 $        
              and from the first relation, we have $ A^{*}_{1} A_{1} A_{2} A^{*}_{2} = 0 $                  
            which gives \be \label{lem9.1} P_{1} P^{{\prime}}_{2} = 0 \ee          
            
            From the second relation, we have $ B^{*}_{1} B_{1} B_{2} B^{*}_{2} = 0 ,$                               
            hence \be \label{lem9.2} Q_{1} Q^{{\prime}}_{2} = 0 \ee          

            Similarly,  the third and fourth relations imply  \be \label{lem9.3}  R_{1} R^{{\prime}}_{2} = 0 \ee
                                  
             and \be \label{lem9.4} ( 1- ( P_{1} + Q_{1}+ R_{1} ) ) ( 1- ( P^{{\prime}}_{2} + Q^{{\prime}}_{2} + R^{{\prime}}_{2} ) ) = 0 \ee         
          respectively.
             
       Now applying the same method to the relations $ A_{1} C_{2} = B_{1} D_{2} = C_{1}A_{2} = D_{1}B_{2} = 0 $ obtained from Lemma \ref{Lemma 8},        
    we obtain \be \label{lem9.5}  P_{1} R^{{\prime}}_{2} = 0 \ee
    
             \be \label{lem9.6} Q_{1} ( 1 - ( P^{{\prime}}_{2} + Q^{{\prime}}_{2} + R^{{\prime}}_{2} ) ) = 0 \ee 
              
             \be \label{lem9.7} R_{1} P^{{\prime}}_{2} = 0 \ee
              
             \be \label{lem9.8} ( 1 - (  P_{1} + Q_{1}+ R_{1} ) ) Q^{{\prime}}_{2} = 0 \ee 
             
      From (\ref{lem9.4}), we get :
      
     $ 1 - ( P^{{\prime}}_{2} + Q^{{\prime}}_{2} + R^{{\prime}}_{2} ) - P_{1} + P_{1} ( P^{{\prime}}_{2} + Q^{{\prime}}_{2} + R^{{\prime}}_{2} ) - Q_{1} + Q_{1} ( P^{{\prime}}_{2} + Q^{{\prime}}_{2} + R^{{\prime}}_{2} ) - R_{1} + R_{1} ( P^{{\prime}}_{2} + Q^{{\prime}}_{2} + R^{{\prime}}_{2} ) = 0 $
      
      Hence, $ 1 - ( P^{{\prime}}_{2} + Q^{{\prime}}_{2} + R^{{\prime}}_{2} ) - P_{1} + P_{1} ( P^{{\prime}}_{2} + Q^{{\prime}}_{2} + R^{{\prime}}_{2} ) - Q_{1} ( 1 - ( P^{{\prime}}_{2} + Q^{{\prime}}_{2} + R^{{\prime}}_{2} ) ) - R_{1} + R_{1} ( P^{{\prime}}_{2} + Q^{{\prime}}_{2} + R^{{\prime}}_{2} ) = 0 $
      
  Applying (\ref{lem9.6}), we have  $ 1 - ( P^{{\prime}}_{2} + Q^{{\prime}}_{2} + R^{{\prime}}_{2} ) - P_{1} + P_{1} ( P^{{\prime}}_{2} + Q^{{\prime}}_{2} + R^{{\prime}}_{2} )  - R_{1} + R_{1} ( P^{{\prime}}_{2} + Q^{{\prime}}_{2} + R^{{\prime}}_{2} ) = 0 $
  
  Now,using (\ref{lem9.2}), we write this as :
  
 $ - ( 1 - ( P_{1} + Q_{1} + R_{1} ) ) Q^{{\prime}}_{2} + 1 - P^{{\prime}}_{2} - R^{{\prime}}_{2} - P_{1} + P_{1} P^{{\prime}}_{2} + P_{1} R^{{\prime}}_{2} - R_{1} + R_{1} P^{{\prime}}_{2} + R_{1} R^{{\prime}}_{2} = 0 $
   
  Now using (\ref{lem9.1}), (\ref{lem9.5}), (\ref{lem9.7}), (\ref{lem9.3}), (\ref{lem9.8}), we obtain 
  
  $ 1 - P^{{\prime}}_{2} - R^{{\prime}}_{2} - P_{1} - R_{1} = 0 $  
  So, we have, $ P_{1} + R_{1} = 1 - ( P^{{\prime}}_{2} + R^{{\prime}}_{2} ) $ \qed

   \brmrk                     
    
     1.From Lemma \ref{Lemma 9} and the fact that $ P_{i} = P^{{\prime}}_{i}, Q_{i} = Q^{{\prime}}_{i}, R_{i} = R^{{\prime}}_{i}, i = 1,2 $ , we have 
             \be \label{1}  P_{1} + R_{1} = 1 - ( P^{{\prime}}_{2} + R^{{\prime}}_{2} ) \ee 
               
              \be \label{2}  P_{1} + R_{1} = 1 - ( P_{2} + R_{2} ) \ee
                 
              \be \label{3}  P^{{\prime}}_{1} + R^{{\prime}}_{1} = 1 - ( P_{2} + R_{2} ) \ee
                 
              \be \label{4}  P^{{\prime}}_{1} + R^{{\prime}}_{1} = 1 - ( P^{{\prime}}_{2} + R^{{\prime}}_{2} )  \ee

            2.  From the above results, we observe that if $ \clq $ is imbedded in $ B( H ) $ for some  Hilbert space H,                    then H breaks up into two orthogonal complements , the  first being the range of $ P_{1} $ and $                         R_{1} $ and the other being the  range of $ Q_{1} $ and $ S_{1} $ .
    
            Let $ p = P^{{\prime}}_1 + R^{{\prime}}_{1} $
            
            Then $p$ is also equal to $ P_{1} + R_{1} = Q^{{\prime}}_{2} + S^{{\prime}}_{2} = Q_{2} + S_{2}$     
            and $ p^{ \bot } = Q^{{\prime}}_{1} + S^{{\prime}}_{1} =  P_{2} + R_{2} = P^{{\prime}}_{2} + R^{{\prime}}_{2} = Q_{1} + S_{1} $ .
   \ermrk

    \blmma :
      \label{Lemma 10}

   $ A_{1} B_{2} - B_{2} A_{1} = 0 = A_{2} B_{1} - { \overline{ \lambda }}^{2} B_{1} A_{2} $   
   
  $ A_{1} D_{2} - { \lambda }^2 D_{2} A_{1} = 0 = A_{2} D_{1} - D_{1} A_{2}  $
    
  $ C_{1} B_{2} - { \lambda }^{2} B_{2} C_{1} = 0 =  B_{1} C_{2} -  C_{2} B_{1} $
   
  $ C_{1} D_{2} - D_{2} C_{1} = 0 =  D_{1} C_{2} - { \lambda }^{2} C_{2} D_{1}  $
  
   \elmma

  {\it Proof :}\\ From Lemma \ref{Lemma 4}, we have  $ A_{1} B_{2} + \overline{ \lambda } B_{1} A_{2} = \lambda A_{2} B_{1} + B_{2} A_{1}. $ 
          So, $ A_{1} B_{2} - B_{2} A_{1} = \lambda (  A_{2} B_{1} - { \overline{ \lambda }}^{2} B_{1} A_{2} ) $        
          Now, $ Ran ( A_{1} B_{2} - B_{2} A_{1} )          
               \subseteq Ran ( A_{1} ) + Ran ( B_{2} )               
               = Ran ( A_{1} A^{*}_{1} ) + Ran ( B_{2} B^{*}_{2} )               
               =  Ran ( P^{{\prime}}_{1} ) + Ran ( Q^{{\prime}}_{2} )               
               \subseteq Ran ( p ) $
               
           On the other hand, $ Ran (  A_{2} B_{1} - { \overline{ \lambda }}^{2} B_{1} A_{2} )           
                               \subseteq Ran ( A_{2} ) + Ran ( B_{1} )                               
                               = Ran ( P^{{\prime}}_{2} ) + Ran ( Q^{{\prime}}_{1} )                               
                               \subseteq Ran ( p^{ \bot } ) $
                               
        So, $  A_{1} B_{2} - B_{2} A_{1} = 0 = A_{2} B_{1} - { \overline{ \lambda }}^{2} B_{1} A_{2} $        
        Similarly, the other three relations can be proved. \qed \\

               Let us now consider  a $C^*$ algebra $\clb$, which has
eight direct summands, four of which are isomorphic with the
commutative algebra $C(\IT^2)$, and the other four are irrational
rotation algebras. More precisely, we take $$\clb=\oplus_{k=1}^8
C^*(U_{k1},U_{k2}), $$   where for odd $k$, 
$U_{k1},U_{k2}$ are the two commuting unitary generators of
$C(\IT^2)$, and  for even $k$, $U_{k1}U_{k2}={\rm exp}(4 \pi i
\theta)U_{k2}U_{k1}$, i.e. they generate $\cla_{2 \theta}$. We set the folowing: $$ \tilde{A_1}:=U_{11}+U_{41},~~\tilde{B_1}:=U_{52}+U_{61},~~\tilde{C_1}:=U_{21}+U_{31},~~\tilde{D_1}:=U_{71}+U_{81},$$
$$\tilde{A_2}:=U_{62}+U_{72},~~\tilde{B_2}:=U_{12}+U_{22},~~\tilde{C_2}:=U_{51}+U_{82},~~\tilde{D_2}:=U_{32}+U_{42}.$$ Denote by $\tilde{M}$ the $4 \times 4$ $\clb$-valued matrix given by $$ \tilde{M}=\left( \begin {array} {cccc}
       \tilde{A_{1}} &  \tilde{A_{2}} & {\tilde{C_{1}}}^* &  {\tilde{C_{2}}}^* \\ \tilde{B_{1}} & \tilde{B_{2}} & {\tilde{D_{1}}}^* & {\tilde{D_{2}}}^* \\ \tilde{C_{1}} & \tilde{C_{2}} & {\tilde{A_{1}}}^* & {\tilde{A_{2}}}^* \\ \tilde{D_{1}} & \tilde{D_{2}} & {\tilde{B_{1}}}^* & {\tilde{B_{2}}}^* \end {array} \right ).$$  We have the following: 
\blmma
\label{bhopf}
(i) The $\ast$-subalgebra generated by the elemements $\tilde{A_i},\tilde{B_i},\tilde{C_i},\tilde{D_i},i=1,2$ is dense in $\clb$;\\
(ii) There is a unique compact (matrix) quantum group structure on $\clb$, where the corresponding coproduct $\Delta_0$, counit $\epsilon_0$ and antipode $\kappa_0$ (say) are  given on the above generating elements by $$ \Delta_0({\tilde{M}}_{ij})=\sum_{k=1}^4 {\tilde{M}}_{ik} \ot {\tilde{M}}_{kj},$$ $$ \kappa_0({\tilde{M}}_{ij})={\tilde{M}_{ji}}^*,~~~\epsilon_0({\tilde{M}}_{ij})=\delta_{ij}.$$
\elmma
The proof can be given by routine verification and hence is omitted. \vspace{2mm}\\
 
Moreover, we have an action of $\clb$ on $\cla_\theta$, as given by the following lemma.
\blmma
\label{baction}
There is a smooth isometric 
 action of $\clb$ on $\cla_\theta$, which is 
given by the following : $$ \alpha_0(U)=U \ot (U_{11}+U_{41})+V
\ot (U_{52}+U_{61})+U^{-1} \ot (U_{21}+U_{31})+V^{-1} \ot
(U_{71}+U_{81}),$$ $$ \alpha_0(V)=U \ot (U_{62}+U_{72})+V \ot
(U_{12}+U_{22})+U^{-1} \ot (U_{51}+U_{82})+V^{-1} \ot
(U_{32}+U_{42}).$$
\elmma
{\it Proof :}\\
It is straightforward to verify that the above indeed defines a smooth action of the quantum group $\clb$ on $\cla_\theta$. To complete the proof, we need to show that $\alpha_0$ keeps the eigenspaces of $\cll$ invariant. For this, we observe that, since $U_{ij}U_{kl}=0$ if $i \neq k$, we have  $$ \alpha_0(U^m)=U^m \ot (U_{11}+U_{41})^m+V^m
\ot (U_{52}+U_{61})^m+U^{-m} \ot (U_{21}+U_{31})^m+V^{-m} \ot
(U_{71}+U_{81})^m,$$ $$ \alpha_0(V^n)=U^n \ot (U_{62}+U_{72})^n+V^n \ot
(U_{12}+U_{22})^n+U^{-n} \ot (U_{51}+U_{82})^n+V^{-n} \ot
(U_{32}+U_{42})^n.$$ From this, it is clear that in the expression of $\alpha_0(U^m) \alpha_0(V^n)$, only coefficients of $U^iV^j$ survive, where $(i,j)$ is one the following:\\ $(m,n), (m,-n), (-m,n), (-m,-n) , (n,m) , (n,-m) ,(-n,m) , (-n,-m).$ This completes the proof of the action being isometric. \qed \\

  Now we are in a position to describe $\clq=QISO(\cla_\theta)$ explicitly.              
          
       \bthm
             $ \clq= QISO ( \cla_{ \theta } )$ is isomorhic (as a quantum group) with $\clb = C ( \IT^{2} ) \oplus \cla_{ 2 \theta } \oplus  C ( \IT^{2} ) \oplus \cla_{ 2 \theta } \oplus C ( \IT^{2} ) \oplus \cla_{ 2 \theta } \oplus  C ( \IT^{2} ) \oplus \cla_{ 2 \theta } $, with the coproduct described before.  
         \ethm   
            
     {\it Proof :}\\ Define $ \phi : \clb \rightarrow \clq $ by 
     
       $$ \phi ( U_{11} ) = A_{1}P^{\prime}_{1}Q^{\prime}_{2} , ~~\phi ( U_{12} ) = B_{2}P^{\prime}_{1}Q^{\prime}_{2} ,~~  \phi ( U_{21} ) = C_{1}{P^{\prime}_{1}}^ \bot Q^{\prime}_{2} ,~~ \phi ( U_{22} ) = B_{2}{P^{\prime}_{1}}^ \bot Q^{\prime}_{2} ,$$
    $$ \phi ( U_{31} ) = C_{1}{P^{\prime}_{1}}^ \bot {Q^{\prime}_{2}}^ \bot ,~~ \phi ( U_{32} ) = D_{2}{P^{\prime}_{1}}^ \bot {Q^{\prime}_{2}}^ \bot,~~\phi ( U_{41} ) = A_{1}P^{\prime}_{1}{Q^{\prime}_{2}}^ \bot ,~~ \phi ( U_{42} ) = D_{2}P^{\prime}_{1}{Q^{\prime}_{2}}^ \bot ,~~   ,$$
    $$   \phi ( U_{51} ) = C_{2}{P^{\prime}_{2}}^ \bot Q^{\prime}_{1} ,~~ \phi ( U_{52} ) = B_{1}{P^{\prime}_{2}}^ \bot Q^{\prime}_{1} ,~~ \phi ( U_{61} ) = B_{1}P^{\prime}_{2}Q^{\prime}_{1} , ~~\phi ( U_{62} ) = A_{2}P^{\prime}_{2}Q^{\prime}_{1},   $$
    $$  \phi ( U_{71} ) = D_{1}P^{\prime}_{2}{Q^{\prime}_{1}}^ \bot,~~\phi ( U_{72} ) = A_{2}P^{\prime}_{2}{Q^{\prime}_{1}}^ \bot ,~~  \phi ( U_{81} ) = D_{1}{P^{\prime}_{2}}^ \bot {Q^{\prime}_{1}}^ \bot,~~\phi ( U_{82} ) = C_{2}{P^{\prime}_{2}}^ \bot {Q^{\prime}_{1}}^ \bot.$$

    We  show that $ \phi $ is well defined and indeed gives a $\ast$-homomorphism. 
    Using the facts that $A_1,B_2$ are commuting normal partial isometries, we have, 
   \bean \lefteqn{  A_{1} P^{\prime}_{1} Q^{\prime}_{2} B_{2} P^{\prime}_{1} Q^{\prime}_{2}}\\     
            &=&  A_{1} A_{1} A^{*}_{1} B_{2} B^{*}_{2} B_{2} A_{1} A^{*}_{1} B_{2} B^{*}_{2}            
            = A_{1} A^{*}_{1} A_{1} B_{2} A_{1} A^{*}_{1} B_{2} B^{*}_{2}\\            
            &=&  A_{1} B_{2} A_{1} A^{*}_{1} B_{2} B^{*}_{2}
            B_{2}{P^{\prime}}_{1} Q^{\prime}_{2} A_{1}{P^{\prime}}_{1} Q^{\prime}_{2}         = B_{2} A_{1} A^{*}_{1} B_{2} B^{*}_{2} A_{1} A_{1} A^{*}_{1} B_{2} B^{*}_{2}\\                                              &=&  A_{1} B_{2} A^{*}_{1} B_{2} B^{*}_{2} A_{1} A^{*}_{1} A_{1} B_{2} B^{*}_{2}                         
                         = A_{1} B_{2} A^{*}_{1} B_{2} B^{*}_{2} A_{1}  B_{2} B^{*}_{2}\\                         & =& A_{1} B_{2} A^{*}_{1} B_{2} A_{1} B^{*}_{2} B_{2} B^{*}_{2}                                               = A_{1} B_{2} A^{*}_{1} A_{1} B_{2} B^{*}_{2}\\ &=& A_{1} B_{2} A_{1} A^{*}_{1} B_{2} B^{*}_{2}.\eean

        So, $ \phi(U_{11})=A_{1} P^{\prime}_{1} Q^{\prime}_{2} $ and $ \phi(U_{12})=B_{2} P^{\prime}_{1} Q^{\prime}_{2} $ commute and they are clearly unitaries when viewed as operators on the range of $P^\prime_1 Q^\prime_2$, which proves that there exists a unique $C^*$-homomorphism from $C(\IT^2)\cong C^*(U_{11},U_{12})$ to $\clq$ which sends $U_{11}$ and $U_{12}$ to $ A_{1} P^{\prime}_{1} Q^{\prime}_{2}$ and $B_{2} P^{\prime}_{1} Q^{\prime}_{2} $ respectively.  
        
   Again, using the facts that $ C_{1} $ and $ B_{2} $ are normal partial isometries satisfying the relation $ B_{2} C_{1} = \frac{1}{{ \lambda }^2} C_{1} B_{2} $, we have,
  $    \phi ( U_{22} ) \phi ( U_{21} ) $  \\
  $ = B_{2} { P^{\prime}_{1} }^{\bot} Q^{\prime}_{2} $ \\
  $ = \frac{1}{\lambda^2} C_{1} {P^{\prime}_{1}}^{\bot} Q^{\prime}_{2} B_{2} {P^{\prime}_{1}}^{\bot} Q^{\prime}_{2} $ \\
  $ = \frac{1}{{\lambda}^2}  \phi ( U_{21} ) \phi ( U_{22} ). $
   
   i.e, $ \phi( U_{21} ) \phi( U_{22} ) = {\lambda}^2 \phi( U_{22} ) \phi( U_{21} ) $ and they are clearly unitaries on the range of $ P^{\prime}_{1} {Q^{\prime}_{2}}^{\bot} $ which proves that there exists a unique $C^*$-homomorphism from $ \cla_{2 \theta} \cong  C^*(U_{21},U_{22})$ to $\clq$ which sends $U_{21}$ and $U_{22}$ to $ C_{1} {{P^{\prime}}_{1}}^{\bot} Q^{\prime}_{2}$ and $B_{2} {{P^{\prime}}_{1}}^{\bot} Q^{\prime}_{2} $ respectively.

    The other cases can be worked out similarly and thus it is shown that $ \phi $ defines a $C^{\ast} $ homomorphism from $ \clb $ to $ \clq $ and moreover,             
   it is easy to see that $\phi(\tilde{M}_{ij})=M_{ij}$, and thus $\phi$ is a morphism of quantum group, and it clearly satisfies $({\rm id} \ot \phi) \circ \alpha_0=\alpha$. By universality of the quantum isometry group $\clq$, this completes the proof that $\clq \cong \clb$ as compact quantum groups.
                         \qed

 \brmrk
   In particular, we note that if $ \theta $ is taken to be $ 1/2 $, then we have a commutative compact quantum group as the quantum isometry group of a noncommutative $ C^{\ast} $ algebra. 
                                                 \ermrk

 We conclude this section with an identification of the `quantum double torus' discovered and studied by Hajac and Masuda (\cite{hajac}) with an interesting quantum subgroup of $QISO(\cla_\theta)$. Consider the $C^*$-ideal $\cli$ of $\clq$ generated by $\tilde{C_i}, \tilde{D_i},$ $i=1,2$. It is easy to verify that $\cli=C^*(U_{ik}, i=2,3,4,5,7,8; k=1,2)$, hence $\clq/\cli \cong C^*(U_{1k},U_{6k},k=1,2)$. Moreover, $\cli$ is in fact a Hopf ideal, i.e. $\clq/\cli$ is a quantum subgroup of $\clq$. Denoting by $A_0,B_0,C_0, D_0$ the elements $U_{11}, U_{61}, U_{62},U_{12}$ respectively, we can describe the structure of $\clq/\cli$ as follows: 
 
 \bthm
 Consider the  $C^*$ algebra
  $\clq^{\rm hol}=C(\IT^2) \oplus \cla_{2\theta}$, given by  the generators $A_0,B_0,C_0,D_0$  ( where $A_0,D_0$ correspond to $C(\IT^2)$ and $B_0,C_0$ correspond to $\cla_{2 \theta}$), with the following coproduct: 
  $$ \Delta_h(A_0)=A_0 \ot A_0+C_0 \ot B_0, ~~\Delta_h(B_0)=B_0 \ot A_0+D_0 \ot B_0,$$
  $$ \Delta_h(C_0)=A_0 \ot C_0+C_0 \ot D_0,~~\Delta_h(D_0)=B_0 \ot C_0+D_0 \ot D_0.$$
  Then $(\clq^{\rm hol},\Delta_h)$ is a compact quantum group isomorphic with $\clq/\cli$.  It has an action $\beta_0$ on $\cla_\theta$ given by
  $$ \beta_0(U)=U \ot A_0+V \ot B_0,~~\beta_0(V)=U \ot C_0 + V \ot D_0.$$ Moreover, $\clq^{\rm hol}$ is universal among the compact quantum groups acting `holomorphically' on $\cla_\theta$ in the following sense: whenever a compact quantum group $(\cls, \Delta)$ has a smooth isometric action $\gamma $ on $\cla_\theta$ satisfying the additional condition that $\gamma$ leaves the subalgebra generated by $\{ U^m V^n,~m,n \geq 0 \}$ invariant, then there is a unique morphism from $\clq^{\rm hol}$ to $\cls$ which intertwins the respective actions.
 \ethm        
 {\it Proof:}\\ We need to prove only the universality of $\clq^{\rm hol}$. Indeed, it follows from the universality of $\clq$ in the category of smooth isometric actions that there is a unique morphism $\phi$ (say) from $\clq$ to $\cls$ such that $\gamma=({\rm id} \ot \phi) \circ \alpha_0$. Writing this relation on $U$, $V$ and noting that by assumption on $\gamma$, the coefficients of $U^{-1},V^{-1}$ in the expression of $\gamma(U), \gamma(V)$ are $0$, it is immediate that $\phi(\tilde{C_i})=\phi(\tilde{D_i})=0$ for $i=1,2$, i.e. $\phi(\cli)=\{0 \}$. Thus, $\phi$ induces a morphism $\tilde{\phi}$ (say) from $\clq/\cli$ to $\cls$ satisfying $\gamma=({\rm id} \ot \tilde{\phi}) \circ \beta_0$. \qed\\

   \section{Quantum isometry group of deformed spectral triples}
    In this section, we give a general scheme for comuputing quantum isometry groups by proving that quantum isometry group of a deformed noncommutative manifold coincides with (under reasonable assumptions) a similar deformation of the quantum isometry  group of the original manifold. To make this precise, we introduce a  few notation and terminologies. 
    
 We begin with some generalities on compact quantum groups.  Given a compact quantum group $(\clg,\Delta)$, recall that the the dense unital $\ast$-subalgebra $\clg_0$ of $\clg$ generated by the matrix coefficients of the irreducible unitary representations has a canonical Hopf $\ast$-algebra structure.  Moreover,  given an action $\gamma : \clb \raro \clb \ot \clg$ of the compact quantum group $(\clg, \Delta)$  on a unital $C^*$-algebra $\clb$, it is known that one can find a dense, unital  $\ast$-subalgebra $\clb_0$ of $\clb$ on which the action becomes an action by the Hopf $\ast$-algebra $\clg_0$( see, for example,  \cite{wang} , \cite{podles} ). We shall use the Sweedler convention of abbreviating $\gamma(b) \in \clb_0 \ot_{\rm alg} \clg_0$ by $b_{(1)} \ot b_{(2)}$, for $b \in \clb_0$. This applies in particular to the canonical action of the quantum group $\clg$ on itself, by taking $\gamma=\Delta$. Moreover, for a linear functional $f$ on $\clg_0$ and an element $c \in \clg_0$ we shall define the `convolution' maps $f \diamond c :=(f \otimes {\rm id} ) \Delta ( c )$ and $c \diamond f := ({\rm id} \otimes f) \Delta ( c)$. We also define convolution of two functionals $f$ and $g$ by $(f \diamond g)(c)=(f \ot g)(\Delta(c))$.

    Let us now consider a $C^*$ dynamical system $(\cla, \IT^n, \beta)$  where    $\beta$ is an action of $\IT^n$, and assume that there exists a spectral triple $(\cla^\infty, \clh, D)$ on the smooth subalgebra $\cla^\infty$ w.r.t. the action of $\IT^n$,   such that the spectral triple satisfies all the assumptions of \cite{goswami} for ensuring the existence of the quantum isometry group. Let $\clq \equiv QISO(\cla)$ denote the quantum isometry group of the spectral triple $(\cla^\infty, \clh,D)$, with $\cll$ denoting the corresponding Laplacian as in \cite{goswami}. Let $\cla_0$ be the $\ast$-algebra generated by complex linear (algebraic, not closed) span $\cla^\infty_0$ of the eigenvectors of $\cll$ which has a countable discrete set of eigenvalues each with finite multiplicities, by assumptions in \cite{goswami}, and it is assumed, as in \cite{goswami}, that $\cla^\infty_0$ is a subset of $\cla^\infty$ and is norm-dense in $\cla$.  Moreover, we make the following assumptions :\\
     (i) $\cla_0$ is dense in $\cla^\infty$ w.r.t. the Frechet topology coming from the action of $\IT^n$.\\ ( ii ) $\bigcap_{n \geq 1 } {\rm Dom}  ({\cll}^n)=\cla^\infty.$\\
     (iii) $\cll$ commutes with the $\IT^n$-action $\beta$, hence $C(\IT^n)$ can be identified as a quantum subgroup of $\clq$.\\
     
      Let $\pi$ denote the surjective map from $\clq$ to its quantum subgroup $C(\IT^n)$, which is a morphism of compact quantum groups. We denote by $\alpha : \cla \raro \cla \ot \clq$ the action of $\clq=QISO(\cla)$ on $\cla$, and note that on $\cla_0$, this action is algebraic, i.e. it is an action of the Hopf $\ast$-algebra $\clq_0$ consisting of matrix elements of finite dimensional unitary representations of $\clq$. We have $ ( id \otimes \pi ) \circ \alpha=\beta$.

  We shall abbreviate  $  e^{ 2 \pi iu}$ by $e(u)$ ($u \in \IR^n$), and shall denote by 
  $ \eta $  the canonical homomorphism from $\IR^n$ to $\IT^n$ given by $ \eta (x_1,x_2,......,x_n ) =(e (x_1),e (x_2),.....e (x_n) )$. For $u \in \IR^n$, $\alpha_u$ will denote the $\IR^n$-action on $\cla$ given by $\alpha_u(a):=({\rm id} \ot \Omega(u))(\alpha(a))$, where    
 $ \Omega ( u ) := {\rm ev}_{\eta ((u))} \circ \pi$, for $u \in R^n$ (${\rm ev}_x$ being the state on $C(\IT^n)$ obtained by evaluation of a function at the point $x \in \IT^n$). 
 
  Let us now briefly recall Rieffel's formulation of deformation quantization (see, e.g. \cite{rieffel}). Let $J$ denote a skew symmetric $n \times n$ matrix with real entries. We define a `deformed' or `twisted' multiplication $\times_J : \cla^\infty \times \cla^\infty \raro \cla^\infty$ given by $$ a \times_J b:=\int \int \alpha_{Ju}(a) \alpha_v(b) e(u.v) du dv,$$ where $u.v$ denotes the standard (Euclidean) inner product on $\IR^n$ and the integral makes sense as an oscillatory integral, described in details in \cite{rieffel} and the references therein. This defines an associative algebra structure on $\cla^\infty$, with the $\ast$ of $\cla$ being an involution w.r.t. the new product $\times_J$ as well, and one can also get a $C^*$-algebra, denoted by $\cla_J$, by completing $\cla^\infty$ in a suitable norm denoted by ${ \left\| \right\|}_{J} $ (see \cite{rieffel}) which is a $C^*$-norm w.r.t. the product $\times_J$. We shall denote by $\cla^\infty_J$ the vector space $\cla^\infty$ equipped with the $\ast$-algebra structure coming from $\times_J$. One has a natural Frechet topology on $\cla^\infty_J$, given by a family of seminorms $ \{ {\left\| \right\|}_{n,J} \} $ where $ \{ {\left\| a  \right\|}_{n,J} \} = \sum_{\left| \mu \right| \leq n} ( \left| \mu \right|!)^{-1} {\left\| \alpha_{X^{\mu}}( a ) \right\|}_{J} , $ ( $ \alpha_{X^{\mu}} $ as in \cite{rieffel}), in which $\cla^\infty_J$ is complete. Moreover, it follows from the estimates ( Proposition 4.10 , page 35 ) in \cite{rieffel} that $\cla^\infty=\cla^\infty_J$ as topological spaces, i.e. they coincide as sets and the corresponding Frechet topologies are also equivalent. In view of this, we shall denote this space simply by $\cla^\infty$, unless one needs to consider it as Frechet algebra, in which case the suffix $J$ will be used. 
  
    Assume furthermore that for each skew-symmetric matrix $J$, there exists a spectral triple on  $\cla^\infty_J$ satisfying the assumptions in \cite{goswami} for defining the quantum isometry group $QISO(\cla_J)$, and assume also that  the corresponding Laplacians, say $\cll_J$, coincide with  $\cll$ on $\cla^\infty \subset \cla_J$, so that the quantum isometry group  $QISO ( \cla_J )$ is the universal compact quantum group acting on $\cla_J$, with the action keeping each of the eigenspaces of $\cll$ invariant. Note that the algebraic span of eigenvectors of $\cll_J$ coincides with that of $\cll$, i.e. $\cla^\infty_0$, which is already assumed to be Frechet-dense in $\cla^\infty=\cla^\infty_J$, hence in particular norm-dense in $\cla_J$. We now state and prove a criterion, to be used later, for extending positive maps defined on $\cla_0$. 
    \blmma
    \label{positive}
    Let $\clb$ be another unital $C^*$-algebra equipped with a $\IT^n$-action, so that we can consider the $C^*$-algebras $\clb_J$ for any skew symmetric $n \times n$ matrix $J$. 
    Let $\phi : \cla^\infty \raro \clb^\infty$ be a linear map, satisfying the following :\\
        (a)  $\phi$  is positive w.r.t. the defomed products $\times_J$ on $\cla_0$ and $\clb^\infty$, i.e. $\phi(a^* \times_J a) \geq 0 $ (in $\clb^\infty_J \subset \clb_J$)  for all $a \in \cla_0$, and \\
    (b) $\phi$ extends to a norm-bounded map (say $\phi_0$) from $\cla$ to $\clb$.\\
    Then $\phi$ also have an extension $\phi_J$   as a $\|~\|_J$-bounded positive map from $\cla_J$ to $\clb_J$ satisfying $\| \phi_J \| =\| \phi(1) \|_J$.
    \elmma
    {\it Proof :}\\
    we can view $\phi$ as a map between the Frechet spaces $\cla^\infty$ and $\clb^\infty$,  which is clearly closable, since it is continuous w.r.t. the norm-topologies of $\cla$ and $\clb$, which are weaker than the corresponding Frechet topologies.  By the Closed Graph Theorem, we conclude that $\phi$ is continuous in the Frechet topology. Since  $\cla^\infty=\cla^\infty_J$ and $\clb^\infty=\clb^\infty_J$ as Frechet spaces, consider $\phi$ as a continuous map from $\cla^\infty_J$ to $\clb^\infty_J$, and  it  follows  by the Frechet-continuity of $\times_J$ and $\ast$ and the Frechet-density of $\cla_0$ in $\cla^\infty_J$ that  the positivity (w.r.t. $\times_J$) of the restriction of $\phi$ to $\cla_0 \subset \cla^\infty_J$ is inherited by the extension on $\cla^\infty=\cla^\infty_J$. Indeed,  given $a \in \cla^\infty_J=\cla^\infty$, choose a sequence $a_n \in \cla_0$ such that $a_n \raro a$ in the Frechet topology. We have $\phi(a^* \times_J a)=\lim_n \phi(a_n^* \times_J a_n)$ in the Frechet topology, so in particular, $\phi(a_n^* \times_J a_n) \raro \phi(a^* \times_J a)$ in the norm of $\clb_J$, which implies that $\phi(a^* \times_J a)$ is a positive element of $\clb_J$ since $\phi(a_n^* \times_J a_n)$ is so for each $n$. 
        Next, we note that $\cla^\infty$ is closed under holomorphic functional calculus as a unital $\ast$-subalgebra of $\cla_J$ (the identity of $\cla^\infty_J$ is same as that of $\cla$), so any positive  map defined on $\cla^\infty$ admits a bounded extension (say $\phi_J$) on $\cla_J$, which will still be a positive map, so in particular the norm of $\phi_J$ is same as $\| \phi_J(1)\|$. \qed \\

  We shall also need Rieffel-type deformation of compact quantum groups (due to Rieffel and Wang, see \cite{wang2},  \cite{toral} and references therein), w.r.t. the action by a quantum subgroup isomorphic to $C(\IT^n)$ for some $n$.   Indeed, for each skew symmetric $n \times n$ real matrix $J$, we can consider a $2n$-parameter action on the compact quantum group, and equip the corresponding Rieffel-deformation the structure of a compact quantum group. We will discuss about it in some more details later on.

    For a fixed $J$, we shall work with several multiplications on the vector space ${\cla_0} \ot_{\rm alg} \clq_0$, where $\clq_0$ is the dense Hopf $\ast$-algebra generated by the matrix coefficients of irreducible unitary representations of the quantum isometry group $\clq$. We shall denote the counit and antipode of $\clq_0$ by $\epsilon$ and $\kappa$ respectively. Let us define the following

    $$ x \odot y = \int_{\IR^{4n}} e( -u.v )e( w.s )(\Omega ( -Ju )\diamond x \diamond (\Omega ( Jw ) ) (\Omega ( -v ) \diamond y \diamond \Omega( s )) du dv dw ds ,$$
     where $x,y \in \clq_0$. This is clearly a bilinear map, and will be seen to be an associative multiplication later on. 
    Moreover, we define two bilinear maps $\bullet$ and $\bullet_J$ by setting  $(a \ot x) \bullet (b \ot y):=ab \ot x \odot y$ and $(a \ot x) \bullet_J (b \ot y):=(a \times_J b) \ot (x \odot y)$, for $a,b \in {\cla_0}$, $x,y \in \clq_0$. 
We have $ \Omega(u) \diamond   ( \Omega(v) \diamond c ) = ( \Omega(u) \diamond \Omega(v) ) \diamond c  $.  
   \blmma 
    \label{Lemma1}
       
  The map $\odot$  satisfies $$ \int_{\IR^{2n}} ( \Omega( J u ) \diamond x ) \odot ( \Omega( v )   \diamond y )e(u.v ) du dv = \int_{\IR^{2n}} (x \diamond  ( \Omega( J u ))( y \diamond \Omega( v ) )e( u.v ) du dv ,$$
   for $x,y \in \clq_0$.
   \elmma
    {\it Proof :}\\
 We have \bean \lefteqn{{\rm LHS}}\\
 &=&   \int(  (  \Omega(Ju^{\prime})) \diamond x ) \odot (\Omega( v^{{\prime}} ) \diamond y )e( u.v ) du^{\prime}dv^{\prime}\\
 & =&  \int_{\IR^{2n}}\{ \int_{\IR^{4n}} e(-u.v )e(w.s )(\Omega(-Ju ) )\diamond ((\Omega( Ju^{\prime})) \diamond x ) \diamond (\Omega( Jw )) ( \Omega( -v )) \diamond (\Omega( v^{\prime} ) \diamond y ) \diamond \Omega( s ) ~    dudvdwds \}e(u^{\prime}.v^{\prime})du^{\prime}dv^{\prime} \\
 &  =& \int_{\IR^{2n}} \{ \int_{\IR^{4n}} e( -u.v )e ( w.s ) (\Omega( -Ju )) \diamond ( \Omega( Ju^{{\prime}} )) \diamond x ) \diamond (\Omega( Jw )) (\Omega( -v )) \diamond (\Omega( v^{{\prime}}) \diamond y ) \diamond \Omega( s ) ~du dv dw ds \} e ( u^{{\prime}} v^{{\prime}} ) du^{{\prime}}dv^{{\prime}}\\
  &  =&  \int_{\IR^{6n}} ( (\Omega( J( u^{\prime}-u )) \diamond x ) \diamond \Omega( Jw )) (\Omega( v^{{\prime}}- v )) \diamond y \diamond \Omega( s ) e( u^{{\prime}}.v^{{\prime}}) e(-u.v )e ( w.s ) du dv dw ds du^{{\prime}} dv^{{\prime}}\\
  &  = & \int_{\IR^{2n}} e( w.s )dw ds \{ \int_{\IR^{4n}} e ( u^{{\prime}}.v^{{\prime}} )e ( -u.v ) du dv du^{{\prime}}dv^{{\prime}} ( \Omega( J ( u^{{\prime}} -u )) \diamond x_w )( ( \Omega( v^{{\prime}} - v )) \diamond y_s ) \},\eean
       where $ x_w = x \diamond \Omega( Jw ) , y_s =y \diamond \Omega( s ) $.

  The proof of the lemma will be complete if we show 
  $$ \int_{\IR^{4n}} e ( u^{{\prime}}.v^{{\prime} } )e ( -u.v ) (\Omega( J( u^{{\prime}} -u )) \diamond x_w ) ( \Omega( v^{{\prime} } -v ) \diamond y_s ) du dv du^{{\prime}} dv^{{\prime}}  = x_w.y_s. $$
  
  By changing variable in the above integral, with  $z = u^{{\prime}} - u,  
       t = v{{\prime}} -v$, it becomes 
     
  $ \int_{\IR^{4n}} e( -u.v ) e ( ( u + z ).( v + t ) ) \phi ( z,t ) du dv dz dt $
  
 $ = \int_{\IR^{4n}} \phi ( z,t ) e ( u.t + z.v ) e ( z.t ) du dv dz dt, $
 where  $$ \phi ( z,t ) = (\Omega( J ( z ) ) \diamond x_w )( \Omega( t ) \diamond y_s ) .$$
  By taking 
  $ ( z,t ) = X, ( v,u ) = Y,$
      and 
   $     F ( X ) = \phi ( z,t ) e( z.t )$,  the integral can be written as 
       \bean \lefteqn{ \int\int F ( X ) e ( X.Y )dX dY }\\
                    & = & F ( 0 )   ~({\rm ~by~ Corollary ~1.12 ~of~ \cite{rieffel},~ page ~9})\\
                     &  =& ( \Omega( J ( 0 ) ) \diamond x_w ) ( \Omega( 0 ) \diamond y_s ) \\
                                        & =& x_w.y_s, \eean since 
                                        \bean  \Omega( J(0) )\diamond x_{w}  = ( ev_{\eta(0)} \pi \otimes id ) \Delta ( x_{w} )  =  ( \epsilon_{\IT^n} \circ \pi \otimes id ) \Delta( x_{w} )  =  ( \epsilon \otimes id ) \Delta( x_{w} )  = x_{w} \eean and similarly  $ \Omega(0) \diamond y_{s} =y_s$  (here $ \epsilon_{\IT^{n}} $ denotes the counit of the quantum group $ C( \IT^{n} ) $).
                                        
                      This proves
                  the claim and hence the lemma. 
   \qed
   
    \blmma
   
  \label{ Lemma2}
   We have for $ a \in {\cla_0}$ 
   $$ \alpha ( \alpha _u ( a ) ) = a_{(1)} \otimes ( id \otimes \Omega( u ) ) ( \Delta ( a_{( 2 )} ).$$
   \elmma
    {\it Proof :}\\
   We have \bean \lefteqn{   \alpha _u ( a )}\\
   & =&  ( id \otimes \Omega ( u ) ) \alpha ( a )\\
                     & =&  ( id \otimes \Omega ( u ) ) ( a_{( 1 )} \otimes a_{( 2 )} )\\
                    &= & a_{( 1 )} ( \Omega ( u )) ( a_{( 2 )} ).\eean
                   
    This gives, \bean \lefteqn{ \alpha ( \alpha _{ u } ( a ) )}\\
           &=&  \alpha ( a_{ ( 1 ) } ) \Omega ( u ) ( a_{ ( 2 ) } )\\
       & = & ( id \otimes id \otimes \Omega ( u ) ) ( \alpha ( a_{(1)} \otimes a_{ (2) } ) )\\
       &= & ( id \otimes id \otimes \Omega ( u ) ) ( ( \alpha \otimes id ) \alpha ( a ) )\\
       &= & ( id \otimes id \otimes \Omega ( u ) ) ( ( id \otimes \Delta ) \alpha ( a ) )\\
       &= & a_{(1)} \otimes ( id \otimes \Omega ( u ) ) \Delta ( a_{(2)}).\eean  
        \qed
        \blmma
       \label{Lemma3}
        For $a,b \in {\cla_0}$, we have      
     $$ \alpha ( a \times _{J} b ) = a_{(1)}b_{(1)} \otimes ( \int\int ( a_{(2)} \diamond  u )  ( b_{(2)} \diamond           v  ) e ( u.v ) du dv  ).$$
     \elmma
     {\it Proof :}\\
     
     Using the notations and definitions in page 4-5 of \cite{rieffel}, we note that for any $ f: \IR^{2} \rightarrow \IC $ belonging to  $ \IB( \IR^{2} ) $ and fixed $ x \in E $( where $ E $ is a Banach algebra ), the function $ F( u, v ) = xf(u,v )$ belongs to  $ \IB^{E} ( \IR^{2} ) $ and we have  \bean \lefteqn{ x ~( \int \int f( u,v ) e( u.v ) du dv )}\\
      & =&  x ~( \lim_L \sum_{p \in L} \int \int ( f \phi_{p})( u,v ) e( u.v ) du dv )\\
    &=& \lim_L \sum_{p \in L} \int \int x~( f \phi_{p})( u,v ) e( u.v ) du dv )\\
      &= & \int\int x~f(u,v ) e(u.v) du dv.\eean
      Then,
        \bean
     \lefteqn{\alpha ( a \times _{J} b )}\\
      &=&  \alpha ( \int\int \alpha _{Ju} ( a ) \alpha _{v} ( b ) e( u.v ) du dv )\\
      &=& \alpha ( \int\int a_{(1)} ( \Omega(Ju) )( a_{(2)} ) b_{(1)} ( \Omega(v) )(b_{(2)}) e(u.v) du dv ) \\
      &=& \alpha (a_{(1)}  b_{(1)}  \int\int  ( \Omega(Ju) )( a_{(2)} )  ( \Omega(v) )(b_{(2)}) e(u.v) du dv ) \\
      &=& \alpha( a_{(1)} ) \alpha( b_{(1)} ) \int\int  ( \Omega(Ju) )( a_{(2)} )  ( \Omega(v) )(b_{(2)}) e(u.v) du dv \\
     &=&  \int \int \alpha( a_{(1)} ) \alpha( b_{(1)} ) ( \Omega(Ju) )( a_{(2)} )  ( \Omega(v) )(b_{(2)}) e(u.v) du dv ) \\
     &=& \int\int \alpha ( \alpha _{Ju}(a) )\alpha ( \alpha _{v} ( b ) ) e( u.v ) du dv \\
      &= & \int\int ( a_{(1)} \otimes ( id \otimes \Omega(Ju) ) ( \Delta ( a_{( 2 )} )) ) ( b_{(1)} \otimes ( id \otimes \Omega ( v ) ) ( \Delta ( b_{( 2 )}) )) e ( u.v ) du dv\\& & ( {\rm ~using ~Lemma ~\ref{ Lemma2} )}\\
     &=&  a_{(1)}b_{(1)} \otimes \int\int ( a _{(2)} \diamond   \Omega( Ju ) ) ( b_{(2)} \diamond  \Omega( v ) ) e ( u.v ) du dv .\eean
     \qed
     
 \blmma    
    \label{Lemma4}
    For $a,b \in {\cla_0}$,      
     $$ \alpha ( a ) \bullet_J \alpha ( b ) = a_{(1)}b_{(1)} \otimes \{ \int\int ( \Omega( Ju )  \diamond a_{(2)} ) \odot (  \Omega( v )  \diamond b_{(2)}) e( u.v ) du dv \}.$$
     \elmma
     {\it Proof :}\\
     We have 
     \bean
     \lefteqn{\alpha ( a )\bullet_J \alpha ( b )}\\
         & = & ( a_{(1)} \otimes a_{(2)} ) ( b_{(1)} \otimes b_{(2)} )\\
        &=&  a_{(1)} \times _{J} b_{(1)} \otimes ( a_{(2)} \odot b_{(2)} )\\
        &=&  \int\int \alpha_{Ju} ( a_{(1)} ) \alpha _{v} ( b_{(1)} ) e( u.v ) du dv \otimes ( a_{(2)} \odot b_{(2)} ). \eean
      Let $\epsilon : \clq_0 \raro \ \IC$ be the counit of the compact quantum group $\clq$.
     So we have $ ( id \otimes \epsilon ) \alpha = id .$    This gives, 
      \bean
      \lefteqn{ \alpha ( a ) \bullet_J \alpha ( b )}\\
     &     =&  \int\int ( id \otimes \epsilon ) \alpha ( \alpha_{Ju} ( a_{(1)} ) ) ( id \otimes \epsilon ) \alpha ( \alpha _{v} ( b_{(1)} ) e( u.v ) du dv \otimes ( a_{2} \odot b_{2} )\\
     &=&  \int\int ( id \otimes \epsilon ) ( \alpha ( \alpha_{Ju} ( a_{(1)} ))) ( id \otimes \epsilon ) ( \alpha ( \alpha_{v} ( b_{(1)} ) ) e( u.v ) du dv \otimes ( a_{(2)} \odot b_{(2)} ). \eean
     Note that by Lemma \ref{Lemma3}, 
     $
     \int\int ( id \otimes \epsilon ) ( \alpha ( \alpha_{Ju} ( a_{(1)} ) ) ( id \otimes \epsilon ) ( \alpha ( \alpha _{v} ( b_{(1)} ) ) ) e(u.v ) du dv $
     
   $  = \int\int ( id \otimes \epsilon ) ( a_{(1)(1)} \otimes ( id \otimes \Omega ( Ju) ) ( \Delta (a_{(1)(2)})) (id \otimes \epsilon ) ( b_{(1)(1)} \otimes ( id \otimes \Omega ( v ) ) ( \Delta ( b_{(1)(2)} )) e( u.v ) du dv $

   $  = \int\int ( id \otimes \epsilon ) ( a_{(1)(1)} \otimes ( a_{(1)(2)} \diamond  \Omega( Ju )  ) ( id \otimes \epsilon ) ( b_{(1)(1)} \otimes ( b_{(1)(2)} \diamond  \Omega( v )  ) e( u.v ) du dv $
     
   $  = \int\int a_{(1)(1)} b_{(1)(1)} \epsilon ( ( a_{(1)(2)} \diamond \Omega( Ju )  ) \epsilon ( b_{(1)(2)} \diamond  \Omega( v )  ) e( u.v ) du dv. $

     Using the fact that  $f \diamond \epsilon=\epsilon \diamond f=f$ for any functional on $\clq_0$, one has 
       $ \epsilon ( a_{(1)(2)} \diamond  \Omega( Ju )  )
                =   \Omega( Ju )   ( a_{(1)(2)} )$ and 
    $ \epsilon ( b_{(1)(2)} \diamond  \Omega( v ) )
                        =  \Omega( v )  ( b_{(1)(2)} )  $, from which it follows that  
                 
        \bean \lefteqn{ \alpha ( a )\bullet_J \alpha ( b )}\\
        &=&  a_{(1)(1)}b_{(1)(1)} \int\int \Omega ( Ju ) ( a_{(1)(2)} ) \Omega ( v ) ( b_{(1)(2)} ) e(u.v ) du dv                    \otimes ( a_{(2)} \odot b_{(2)} )\\
       &=&  \int\int ( id \otimes \Omega ( Ju ) \otimes id ) ( a_{(1)(1)} \otimes a_{(1)(2)} \otimes a_{(2)} )            \bullet ( id \otimes \Omega ( v ) \otimes id ) ( b_{(1)(1)} \otimes b_{(1)(2)} \otimes b_{(2)} )\\
       & &  e( u.v ) du dv\\
  &=&  \int\int  ( id \otimes \Omega ( Ju ) \otimes id ) ( a_{(1)} \otimes \Delta ( a_{(2)} ) ) \bullet ( id \otimes ( \Omega ( v ) \otimes id ) ( b_{(1)} \otimes \Delta ( b_{(2)} ) ) e( u.v ) du dv \\
   &   =&  \int\int \{  a_{(1)} \otimes ( \Omega ( Ju ) \otimes id ) \Delta ( a_{2} ) \} \bullet \{  b_{1} \otimes ( \Omega ( v ) \otimes id ) \Delta ( b_{2} ) \} e( u.v ) du dv\\
      &=&  a_{(1)}b_{(1)} \otimes \int\int ( ( \Omega ( Ju ) \otimes id ) \Delta ( a_{2} ) ) \odot ( \Omega ( v ) \otimes id ) ) \Delta ( b_{(2)} ) e( u.v ) du dv\\
    &=& a_{(1)}b_{(1)} \otimes \int\int (  \Omega( Ju )  \diamond a_{(2)} ) \odot (  \Omega( v )  \diamond b_{(2)} ) e( u.v ) du dv,\eean
              where we have used the relation $(\alpha \ot id) \alpha=(id \ot \Delta) \alpha$ to get  $a_{(1)(1)} \otimes a_{(1)(2)} \otimes a_{(2)}
   =a_{(1)} \otimes \Delta ( a_{(2)} )$ 
                            and similarly $ b_{(1)(1)} \otimes b_{(1)(2)} \otimes b_{(2)} = b_{(1)} \otimes \Delta ( b_{(2)} )   .$ \qed

        Combining Lemma \ref{Lemma1}, Lemma \ref{Lemma3} and Lemma \ref{Lemma4} we conclude the following.      
   \blmma           
         \label{Lemma5}
         For $a,b \in {\cla_0},$ we have $\alpha(a) \bullet_J \alpha(b)=\alpha(a \times_J b).$
         \elmma

  We shall now identify $\odot$ with the multiplication of a Rieffel-type deformation of $\clq$. Since $\clq$ has a quantum subgroup isomorphic with $\IT^n$, we can consider the following canonical action $\lambda$ of $\IR^{2n}$ on $\clq$ given by          
      $$  \lambda_{( s,u )} = ( \Omega( -s ) \otimes id ) \Delta ( id \otimes \Omega ( u ) ) \Delta.$$  
          Now, let 
       $ \widetilde J := -J \oplus J $, which is a skew-symmetric ${2n} \times 2n$ real matrix, so one can deform $\clq$ by defining the product of $x$ and $y$ ($x,y \in \clq_0$, say) to be the following: $$ \int\int \lambda_{\widetilde J ( u,w )}( x ) \lambda _{v,s} ( y ) e ( ( u,w ).( v,s ) ) d ( u,w ) d ( v,s ).$$ We claim that this is nothing but $\odot$ introduced before.
         
         \blmma
              
           $ x \odot y = x \times_{\widetilde{J}} y ~ \forall x,y \in Q_{0} $
           \elmma   
          {\it Proof :}\\
           Let us first observe that 
          \bean
          \lefteqn{ \lambda _{ \widetilde J ( u,w ) } ( x )}\\
           &=& ( \Omega ( Ju ) \otimes id ) \Delta ( id \otimes \Omega ( Jw ) ) \Delta ( x )\\
          &=& \Omega( Ju ) \diamond x \diamond \Omega( Jw ), \eean
          and similarly
           $ \lambda_{(v,s)}(y)= \Omega( -v ) \diamond y \diamond \Omega( s ).$ 
                     
     Thus, we have                
     \bean     
          \lefteqn{ x \odot y}\\
          & =&  \int_{\IR^{4n}} (  \Omega( -Ju ) \diamond x \diamond  \Omega( Jw ) ) (  \Omega( -v ) \diamond y \diamond                  \Omega( s ) ) e ( -u.v ) e ( w.s ) du dv dw ds\\
      &=&  \int_{\IR^{4n}}   \Omega( J u^{{\prime}} ) \diamond x \diamond  \Omega( J w ) ) \Omega( -v )               \diamond y \diamond  \Omega( s )) e ( u^{{\prime}}.v ) e ( w.s ) du^{{\prime}} dv dw ds\\
      &=& 
     \int_{\IR^{2n} } \int_{\IR^{2n}} \lambda _{\widetilde J ( u,w )} ( x )\lambda _{( v,s )} ( y ) e ( ( u,w ).( v,s )            ) d( u,w ) d ( v,s ),\eean which proves the claim. \qed \\

     Let us denote by $\clq_{\widetilde{J}}$ the $C^*$ algebra obtained from $\clq$ by the  Rieffel deformation w.r.t. the matrix $\widetilde{J}$ described above. It has been shown in \cite{wang2} that the coproduct $\Delta$ on $\clq_0$ extends to a coproduct for the deformed algebra  as well and  $(\clq_{\widetilde{J}}, \Delta)$ is a compact quantum group.
     
      \blmma           
         \label{Lemma5a}
         
     The Haar state (say $h$) of $\clq$ coincides with the Haar state on $\clq_{\widetilde{J}}$ ( say $ h_{J} $ ) on the common subspace $\clq^{\infty}$, and moreover, $h(a \times_{\widetilde{J}} b)=h(ab)$ for $a,b \in \clq^{\infty}$.
         \elmma
         
      {\it Proof :}\\ From \cite{wang2} ( Remark 3.10(2) ), we have that $h$ = $h_{J}$ on $\clq_{0}$. 
                       By using  $ h ( \Omega ( - s ) \otimes id ) = \Omega ( - s )( id \otimes h ) $ and $ h( id \otimes \Omega( u ) ) = \Omega ( u ) ( h \otimes id ) $, we  have for $a \in \clq_0$ 
 
 \bean \lefteqn { h  ( \lambda_{s,u} ( a ) ) }\\
            & = & \Omega ( - s ) ( id \otimes h ) \Delta ( id \otimes \Omega ( u ) ) \Delta ( a ) \\
            & = & \Omega ( -s ) ( h ( ( id \otimes \Omega ( u ) ) \Delta ( a ) )1 ) \\
            & = & h ( ( id \otimes \Omega ( u ) ) \Delta ( a ) ) \\
            & = & \Omega ( u ) ( h( a ).1 ) \\
            & = & h( a ) .\eean             
     
    Therefore,  $ h \lambda_{s,u} ( b ) = h( b ) ~ \forall b \in \clq_0 $.  
      Now,
      \bean
      \lefteqn{ h( a {\times}_{\widetilde{J}} b )}\\
      & =& \int\int h( \lambda_{\widetilde{J}u} ( a ) \lambda_{v} ( b ) ) e( u . v ) du dv \\
      & =& \int\int h( \lambda_{v} ( \lambda_{\widetilde{J}u - v} ( a ) b ) ) e( u . v ) du dv \\
      &=& \int\int h( \lambda_{t} ( a ) b ) e( s.t ) ds dt ,\eean where $ s = - u, t = \widetilde{J}u - v ,$
   which by Corollary 1.12, \cite{rieffel} equals $ h( \lambda_{0} ( a ) b ) = h( a b ). $  That is, we have proved \be \label{etaseta} \lgl a, b \rgl_J=\lgl a, b \rgl ~~\forall a,b \in \clq_0, \ee
      where $\lgl \cdot, \cdot \rgl_J$ and $\lgl \cdot, \cdot \rgl$ respectively denote the inner products of $L^2(h_J)$ and $L^2(h)$. 
    We now complete the proof of the lemma by extending (\ref{etaseta}) from   $\clq_0$ to $\clq^\infty$, by using the fact that $\clq^\infty$ is a common  subspace of the Hilbert spaces $L^2(h)$ and $L^2(h_J)$ and moreover, $\clq_0$ is  dense in both these Hilbert spaces. In particular, taking $a=1 \in \clq_0$, we have $h=h_J$ on $\clq^\infty$.   
      \qed

    \brmrk
    \label{haarrem}
    
    Lemma \ref{Lemma5a} implies in particular that for every fixed $a_1,a_2 \in \clq_0$, the functional $\clq_0 \ni b \mapsto h(a_{1} \times_{\widetilde{J}} b \times_{\widetilde{J}} a_2)=h((\kappa^2(a_2)\times_{\widetilde{J}}a_1)b)$ extends to a bounded linear functional on $\clq$.
\ermrk

 \blmma
     \label{Lemma5b}
     
   If $h$ is faithful on $ \clq $, then $ h_{J} $ is faithful on $\clq_{\widetilde{J}}$.
   
   \elmma
   
     {\it Proof :}\\ Let $a \geq 0, \in \clq_{\widetilde{J}} $ be such that $ h( a ) = 0$.
     
     Let $e$ be the identity of $ \IT^{2n} $ and $U_{n}$ be a sequnce of neighbourhoods of $e$ shrinking to $e$, $f_{n}$  smooth, positive functions with support contained inside $U_{n}$ such that $ \int f_{n}(z) dz = 1 ~ \forall n$.
     Define $ \lambda_{f_{n}} ( a ) = \int_{\IT^{2n}} \lambda_{z} ( a ) f_{n}( z ) dz $.  It is clear that  $ \lambda_{f_{n}} ( a ) \in \clq^{\infty}$ and is positive in $ \clq_{\widetilde{J}}$. Moreover, using the fact that the map  $ z \mapsto \lambda_{z} ( a ) $ is continuous $ \forall a $, we have $ \lambda_{f_{n}} ( a ) \rightarrow a $ as $ n \rightarrow \infty $.
     Now \bean \lefteqn {h_{J}( \lambda_{f_{n}} ( a ) )} \\
                       & = & \int_{\IT^{2n}} h_{J} ( \lambda_{z} ( a ) ) f_{n} ( z ) dz \\
                       & = & \int_{\IT^{2n}} h_{J} ( a ) f_{n} ( z ) dz \\
                       & = & 0, \eean
           so we have
                      $ h( \lambda_{f_{n}} ( a ) ) = 0 $, since $h$ and $h_J$ coincide on $\clq^\infty$ by Lemma \ref{Lemma5a}.
         
 Now we fix some notation which we are going to use in the rest of the proof. Let $ L^{2}( h ) $ and $ L^{2}( h_{J} ) $ denote the G.N.S spaces of $ \clq $ and $ \clq_{\widetilde{J}} $ respectively with respect to the Haar states. Let $ i $ and $ i_{J} $ be the canonical maps from $ \clq $ and $ \clq_{\widetilde{J}} $ to  $ L^{2}( h ) $ and $ L^{2}( h_{J} ) $ respectively. Also, let $ \Pi_{J} $ denote the G.N.S representation of $ \clq_{J} $.
  Using the facts  $ h ( b^{\ast} \times_{\widetilde{J}} b ) = h ( b^{\ast} b ) ~ \forall b \in \clq^{\infty} $ and $ h = h_{J} $ on $ \clq^{\infty} $ ,  we get $ \left\|i_{J}( b ) \right\|^{2}_{L^{2}( h_{J} )} = \left\|i( b ) \right\|^{2}_{L^{2}( h ) } \forall b \in \clq^{\infty} $. So the map sending $ i ( b ) $ to $ i_{J} ( b ) $ is an isometry from a dense subspace of $ L^{2}( h )$ onto a dense subspace of $ L^{2}( h_{J} ) $, hence it extends to a unitary, say $ \Gamma : L^{2} ( h ) \rightarrow L^{2} ( h_{J} ) $. We also note that the maps $  i $ and $ i_{J} $ agree on $ \clq^{\infty} $.
  
  Now, $ \lambda_{f_{n}} ( a ) = b^{\ast} \times_{\widetilde{J}} b $ for some $ b \in \clq_{\widetilde{J}} $. So $ h( \lambda_{f_{n}} ( a ) ) = 0 $ implies $ \left\|i_{J} ( b ) \right\|^{2}_{L^{2} ( h_{J} ) } = 0 $. Therefore, one has  $ \Pi_{J}( b^{\ast} ) i_{J} ( b ) = 0$, and hence $ i_{J}( b^{\ast} b ) = i_{J} ( \lambda_{f_{n}} ( a ) ) = 0$. It thus follows that $ \Gamma ( i ( \lambda_{f_{n}}( a ) ) ) = 0$, which implies  $ i ( \lambda_{f_{n}}( a ) ) = 0$. But the faithfulness of $h$  means that $ i $ is one one, hence $ \lambda_{f_{n}}( a ) = 0$ for all $n$.
 Thus,   $a=lim_{n \rightarrow \infty} \lambda_{f_{n}} ( a ) =0 $, which proves the faithfulness of $h_J$. \qed

     \bthm
     \label{premain}
   If the Haar state is faithful on $\clq$, then $\alpha : {\cla_0} \raro {\cla_0} \ot \clq_0$ extends to an action of the compact quantum group $\clq_{\widetilde{J}}$ on $\cla_J$, which is isometric, smooth and faithful. 
     \ethm
        {\it Proof :}\\
        We have already seen in Lemma \ref{Lemma5} that $\alpha$ is an algebra homomorphism 
     from $\cla_0$ to $\cla_0 \ot_{\rm alg} \clq_0$ (w.r.t. the deformed products), and it is also $\ast$-homomorphism since it is so for the undeformed case and the involution $\ast$ is the same for the deformed and undeformed algebras.  
     It now suffices to show that $\alpha$ extends to $\cla_J$ as a $C^*$-homomorphism. 
     Let us fix any faithful imbedding $\cla_J \subseteq \clb(\clh_0)$ (where $\clh_0$ is a Hilbert space) and consider the imbedding $\clq_{\widetilde{J}} \subseteq \clb(L^2(h_J))$. By definition, the norm on $\cla_J \ot \clq_{\widetilde{J}}$ is the minimal (injective) $C^*$-norm, so it is equal to the norm inherited from the imbedding $\cla_J \ot_{\rm alg} \clq_{\widetilde{J}} \subseteq \clb(\clh_0 \ot L^2(h_J))$. Let us consider the dense subspace $\cld \subset \clh_0 \ot L^2(h_J)$ consisting of vectors which are finite linear combinations of the form $\sum_i u_i \ot x_i,$ with $u_i \in \clh_0$, $x_i \in \clq_0 \subset L^2(h_J)$. Fix such a vector $\xi=\sum_{i=1}^k u_i \ot x_i$ and consider $\clb:=\cla \ot M_k(\IC)$, with the $\IT^n$-action $\beta  \otimes \it {\rm id}$ on $\clb$. Let $\phi: \cla^\infty \raro \clb^\infty$ be the map given by $$ \phi(a):= \left( \left( ({\rm id} \ot \phi_{(x_i,x_j)})(\alpha(a)) \right) \right)_{ij=1}^k,$$ where $\phi_{(x,y)}(z):=h(x^* \times_{\widetilde{J}} z \times_{\widetilde{J}} y)$ for $x,y,z \in \clq_0$. 
     Note  that the range of $\phi$ is in $\clb^\infty=\cla^\infty \ot M_k(\IC)$ since we have  $\phi_{(x,y)}(\cla^\infty) \subseteq \cla^\infty$ by the Remark 2.16 of \cite{goswami}, using our assumption (ii) that $\bigcap_{n \geq 1}{\rm Dom}(\cll^n)=\cla^\infty$.
     
      Since  $\alpha$ maps $\cla_0$ into $\cla_0 \ot_{\rm alg} \clq_0$ and $h=h_J$ on $\clq_0$, it is easy to see that   for $a \in \cla_0$, $\phi(a^* \times_J a)$ is positive in $\clb_J$. Moreover,     by the Remark \ref{haarrem}, $\phi_{(x_i,x_j)}$ extends to $\clq$ as  a bounded linear  functional, hence $\phi$ extends to a bounded linear (but not necessarily positive) map from $\cla$ to $\clb$.  
      Thus, the hypotheses of  Lemma \ref{positive} are satisfied and we conclude that $\phi$ admits a positive extension, say $\phi_J$, from $\cla_J$ to $\clb_J=\cla_J \ot M_k(\IC)$. Thus, we have for $a \in \cla_0$, \bean \lefteqn{\sum_{i,j=1}^k \lgl u_i, \phi(a^* \times_J a)u_j \rgl}\\
      & \leq & \| a \|_J^2 \sum_{ij} \lgl u_i, \phi(1)u_j \rgl 
      = \| a\|_J^2 \sum_{ij} \lgl u_i,u_j \rgl h(x_i^* \times_{\widetilde{J}} x_j)\\
      &=& \| a \|_J^2 \sum_{ij} \lgl u_i \ot x_i, u_j \ot x_j \rgl =\|a\|_J^2 \| \sum_{i=1}^k u_i \ot x_i \|^2.\eean This implies $$ \| \alpha(a) \xi \|^2=\lgl \xi, \alpha(a^* \times_J a) \xi \rgl \leq \| a \|^2_J \| \xi \|^2$$ for all $\xi \in \cld$ and $a \in \cla_0$,  hence $\alpha$ admits a bounded extension which is clearly a $C^*$-homomorphism.    
                 \qed \\
     
  Let $ \clc_{J} $ be the category of compact quantum groups acting isometrically on $ \cla_{J} $ with objects being the pair $ ( \cls,\alpha_{\cls} ) $ where the compact quantum group $ \cls $ acts isometrically on $ \cla_{J} $ by the action $ \alpha_{\cls} . $ If the action is understood, we may simply write $ ( \cls,\alpha_{\cls} ) $ as $ \cls . $  For any two compact quantum groups $ \cls_1 $ and $ \cls_2 $ in $\clc_J$, we write $ \cls_1 < \cls_2 $ if   there is a surjective $C^*$ homomorphism $\pi$ from $ \cls_2 $ to $ \cls_1 $ preserving the respective coproducts (i.e. $\cls_1$ is a quantum subgroup of $\cls_2$) and $\pi$ also satisfies $\alpha_{\cls_1}=({\rm id} \ot \pi) \circ \alpha_{\cls_2}.$   
     
      \brmrk
           
          \label{pqrs} 
           It can easily be seen that $ S_1 < S_2 $ means that $ ( S_1 )_J < ( S_2 )_J $ 
        \ermrk
        
           \bthm          
           \label{abcd}
           If the Haar state on $QISO( \cla )$ is faithful, we have the isomorphism of compact quantum groups:
           $$( QISO ( A ) )_{ \widetilde J } \cong  QISO ( A_{J} ). $$ 
           
           \ethm
           
        {\it Proof :}\\  Let $ \clq ( \cla_{J} ) $ is the universal object in $ \clc_{J} $. By Theorem \ref{premain},  we have seen that $ ( \clq ( \cla ) )_{\widetilde J} $ also acts faithfully, smoothly and isometrically on $ \cla_{J} $, which implies, 
           
             $$ ( \clq ( \cla ))_{\widetilde  J} < \clq ( \cla_{J} ) ~~{\rm in}~ \clc_J .$$
           
           So, by Remark \ref{pqrs}, $ ( ( \clq ( \cla )_{\widetilde J})_{- \widetilde J} < ( \clq ( \cla_{J} ) )_{- \widetilde J} $ in $ \clc_{0} $, hence  
           $ \clq ( \cla ) < ( \clq ( \cla_{J} ) )_{- \widetilde J}.$
           
           Replacing $ \cla $ by $ \cla_{-J} $, we have 
           \bean \lefteqn{ \clq ( \cla_{-J} )}\\
           & <&  \clq ( ( \cla_{- J} )_{J} )_{- \widetilde J} ~~({\rm  in}~~  \clc_{-J})\\
            & \cong &  \clq ( \cla )_{- \widetilde J} ~~({\rm  in}~~  \clc_{-J})~
                           \cong ( \clq ( \cla ) )_{\widetilde {-J}}. \eean                 
           Thus,  $ \clq ( {\cal A}_{J} ) < ( \clq ( \cla ) )_{\widetilde J} $ in $ \clc_{J}, $ which implies 
            $ \clq ( \cla_{J} ) \cong ( \clq ( \cla ) )_{\widetilde J} $ in $ \clc_{J}. $   \qed \\ 
   \bxmpl
            
      We recall that $\cla_{\theta}$ is a Rieffel type deformation of $ C( \IT^2 ) $,( See \cite{rieffel},example 10.2, page 69 ) and it can be easily verified that in this case the hypotheses of this section are true. So  Theorem \ref{abcd} can be applied to compute $ QISO ( \cla_{\theta} ) .$ This gives an alternative way to prove the results obtained in subsection 2.3.                    
           
     \exmpl 
     
      \bxmpl
   We can apply our result to the isospectral deformations of compact oriented Riemannian manifolds considered in  \cite{CD}, in particular to the deformations $S^n_\theta$ of classical $n$-sphere, with the spectral triple defined in \cite{CD}. Since we have proved that $QISO(S^n) \cong C(O(n))$, it will follow that $QISO(S^n_\theta) \cong O_\theta(n)$, where $O_\theta(n)$ is the compact quantum group obtained in \cite{CD} as the $\theta$-deformation of $C(O(n))$. 
   \exmpl  
     
     \brmrk
     
          We would like to conclude this article with the following important and interesting open question : Does there exist a connected, compact manifold whose quantum isometry group is non commutative as a $ C^{*} $ algebra ? We have already observed that for $ S^{n}, \IT^{1}, \IT^{2} $, the answer is negative.
          
          \ermrk
          
 {\noindent}   {\bf Acknowldgement:}\\
     The authors thank P. Hajac and S.L. Woronowicz for some stimulating discussion.

\end{document}